# All Mach Number Second Order Semi-Implicit Scheme for the Euler Equations of Gasdynamics

S. Boscarino [*], G. Russo [†] L. Scandurra [‡]

**Abstract**

This paper presents an asymptotic preserving (AP) all Mach number finite volume shock capturing method for the numerical solution of compressible Euler equations of gas dynamics. Both isentropic and full Euler equations are considered. The equations are discretized on a staggered grid. This simplifies flux computation and guarantees a natural central discretization in the low Mach limit, thus dramatically reducing the excessive numerical diffusion of upwind discretizations. Furthermore, second order accuracy in space is automatically guaranteed. For the time discretization we adopt an Semi-IMplicit/EXplicit (S-IMEX) discretization getting an elliptic equation for the pressure in the isentropic case and for the energy in the full Euler equations. Such equations can be solved linearly so that we do not need any iterative solver thus reducing computational cost. Second order in time is obtained by a suitable S-IMEX strategy taken from Boscarino et al. in [6]. Moreover, the CFL stability condition is independent of the Mach number and depends essentially on the fluid velocity. Numerical tests are displayed in one and two dimensions to demonstrate performances of our scheme in both compressible and incompressible regimes.

# 1 Introduction

Numerical methods for the solution of hyperbolic systems of conservation laws has been a very active field of research in the last decades. Several very effective schemes are nowadays treated in textbooks which became a classic on the topic [33, 45, 21]. Because of the hyperbolic nature, all such systems develop waves that propagate at finite speeds. If one wants to accurately compute all the waves in a hyperbolic system, then one has to resolve all the space and time scales that characterize it. Most schemes devoted to the numerical solution of such systems are obtained by explicit time discretization, and the time step has to satisfy a stability condition, known as CFL condition, which states that the time step should be limited by the space step divided by the fastest wave speed (times a constant of order 1). Usually such a restriction is not a problem: because of the hyperbolic nature of the system, if the order of accuracy is the same in space and time, accuracy restriction and stability restrictions are almost the same, and the system is not *stiff*. There are, however, cases in which some of the waves are not particularly relevant and one is not interested in resolving them. Let us consider as a prototype model the classical Euler equations of compressible gas dynamics. In the low Mach number regimes, it may happen that the acoustic waves carry a negligible amount of energy, and one is mainly interested in accurately capturing the motion of the fluid. In such a case the system becomes *stiff*: classical CFL condition on the time step is determined by the acoustic waves which have a negligible influence on the solution, but which deeply affect the efficiency of the method itself.

[*]Department of Mathematics and Computer Science, University of Catania, Catania, 95125, E-mail: boscarino@dmi.unict.it

[†]Department of Mathematics and Computer Science, University of Catania, Catania, 95125, E-mail: russo@dmi.unict.it

[‡]Heinrich-Heine-Universität Düsseldorf, North Rhine-Westphalia, Germany, E-mail: scandurr@uni-duesseldorf.de

Another difficulty arising with standard Godunov-type schemes for low-Mach flows is that the amount of numerical viscosity on the slow waves introduced by upwind-type discretization of the system would heavily degrade the accuracy. An account of the latter effect is analyzed in [20], where the relevance of centering pressure gradients in the limit of small Mach number is emphasized.

In order to overcome the drawback of the stiffness, one has to resort to implicit strategies for time discretization, which avoid the acoustic CFL restriction and allow the use of a much larger time step. Naive implementation of implicit schemes for the solution of the Euler equations presents however two kinds of problems. First, classical upwind discretization (say Godunov methods based on exact or approximate Riemann solvers) are highly nonlinear and very difficult to solve implicitly. Second, the implicit version of classical schemes may introduce an excessive numerical dissipation on the slow wave, resulting in loss of accuracy. Investigation of the effect on fully implicit schemes (and preconditioning techniques adopted to cure the large numerical diffusion) are discussed for example in [47] and in [34], both inspired by an early work by Turkel [46]. In both cases, a modification to the absolute value of the Roe matrix is proposed by a suitable preconditioner that avoids excessive numerical diffusion of upwind-type discretization at very low Mach.

Several techniques have been devised to treat problems in the low Mach number regimes, that alleviate both drawbacks (see for example [31]). However, some of such techniques have been explicitly designed to treat low Mach number regimes, and are based on low Mach number asymptotics ([27], [29]). There are cases in which the Mach number can change of several orders of magnitude. The biggest challenges come from gas dynamic problems in astrophysics, where the range of scales of virtually all parameters vary over several orders of magnitude. An adaptive low Mach number scheme, based on a non conservative formulation, has been developed with the purpose of tackling complex gas dynamics problems in astrophysics (see [38] and references therein). When Mach number is very low the flow does not develop shock discontinuities, and the conservation form of the schemes is not mandatory. When Mach number is not small, then shock discontinuities may form. In such a cases it is necessary to resort to conservative schemes (see for example [34] for other astrophysical applications).

Several other physical systems are affected by drastic changes of the sound speed. Such large variation may be due to geometrical effects, as for example in the case of the nozzle flow [2] or to heterogeneity of the media. Air-water systems, for example, are characterized by density ratio of three orders of magnitude, while the ratio of sound speed is about five. Waves in heterogeneous solid materials may travel at very different speeds, depending on the local stiffness of the medium. The motivation for the construction of effective all Mach number solver is twofold: on one hand it is relevant to be able to accurate simulate waves in heterogeneous materials without small time step restriction suffered by explicit schemes, on the other hand such simulations can be adopted as a tool to validate homogenized models, which at a more macroscopic scale can be described as a homogeneous medium with different mechanical properties. For example, in air-water flows, for a range of values of the void fraction the measured sound speed is lower than both water and air sound speed [17].

Motivated by the above arguments, several researchers have devoted a lot of effort in the development of *all Mach number* solvers for gas dynamics. An attempt in this direction is presented in [32], where the authors adopt a pressure stabilization technique to be able to go beyond the classical CFL restriction. The technique works well for moderate Mach number, but is not specifically designed to deal with very small Mach numbers.

In an impressive sequence of papers and conference proceedings, [15, 11, 13, 14, 12], F.Coquel and collaborators proposed a semi-implicit strategy, coupled with a multi resolution approach, for the numerical solution of hyperbolic systems of conservation laws with well separated wave propagation speeds. In particular, they considered application to fluid mixtures, in which the propagation speed of acoustic waves, often carrying a negligible amount of energy, is much larger than the speed of the material wave traveling at the fluid velocity. The basic framework is set in [15]. The method is first explained in the context of linear hyperbolic systems. The eigenvalues are sorted and it is assumed that there is a clear separation between slow and fast waves. The Jacobian matrix is split into a slow and fast component, using the characteristic decomposition. The flux at cell boundaries is consequently split into a slow and fast term. The fast term is treated implicitly, while the slow one is treated explicitly. The approach is then generalized to the quasilinear case, making use of Roe-type approximation of flux difference. This allows to construct a simple semi-implicit formulation by leaving the Roe matrix of the fast waves at the previous time step, while only the field is computed at the new time step, leading to a linearly implicit scheme. The effectiveness of the approach is further improved by adopting spatial multi resolution: given a multi scale expansion of the numerical solution, the finest scale is maintained locally only where needed, while coarser scales are adaptively adopted in smoother regions, with a great savings in computational time. Different schemes, still adopting implicit-explicit time differentiation to filter out fast waves, are considered in [11], where a sort of arbitrary Lagrangian-Eulerian scheme is constructed: a fractional time step strategy is composed by an implicit Lagrangian step, which filters out acoustic waves, and an explicit Eulerian step, which takes into account the contribution of slow waves. The main application is still on a model for the evolution of gas-oil mixture. In order to simplify the treatment of a general equation of state, a relaxation method is adopted (which of course satisfies the Chen-Levermore-Liu sub-characteristic condition [10]). The problem of developing an adaptive (local) time step strategy is considered in the proceedings [13], and fully exploited in [14]. In [12], the authors further refine the technique, thus producing a positivity preserving, entropic semi-implicit scheme for Euler-like equations. The approach developed by Coquel and collaborators is certainly valuable, although it may be quite involved to be efficiently implemented for more complex, multidimensional situations.

A different approach has been adopted by Munz and collaborators, starting from the low Mach number asymptotic of Kleinerman and Majda [28, 30]. In [35], the authors develop a very effective semi-implicit method which can be viewed as a generalization of a compressible solver to weakly compressible flows. The method is based on the asymptotic behavior of the Euler equation for low Mach number. Two pressures are defined, a thermodynamic one, which is essentially constant in space, and a dynamic one, which accounts for fluid motion. The method is based on a discretization of the system written in primitive variables. The approach, designed for low Mach flow, cannot be directly used when compressive effects are more pronounced. In a subsequent paper [41], Park and Munz extend the method, still using the pressure as basic unknown in place of the energy, but now they adopt a conservative formulation, thus being able to capture shocks when the Mach number is not so small. Several space discretizations as well as time discretization strategies are discussed, which allow to obtain second order accuracy in space and time. In addition, the paper contains a nice overview of other works on low Mach number flow.

In [22] and in [19] the authors explore the construction of an all Mach-number finite volume scheme for the isentropic Euler and Navier-Stokes equations. In both cases, the approach consists in a sort of hyperbolic splitting, obtained by adding and subtracting a gradient-type term to the momentum equation. Such a term is an approximation of the pressure gradient, and is treated implicitly, while the (relatively small) difference with the physical pressure gradient is treated explicitly. The authors show the asymptotic preserving (AP) property of the schemes: when the Mach number approaches zero the schemes become a consistent and stable discretization of the incompressible Euler and Navier-Stokes equations. In a more recent paper, Cordier et al.[16] extend the technique to the full Euler and Navier-Stokes equations. In paper [18] a different approach has been adopted for the construction of asymptotic preserving schemes for the gas dynamics. The authors perform a *gauge* decomposition of the momentum density into a solenoidal and irrotational field. They show that this corresponds to a sort of micro-macro decomposition, in which the macroscopic variable describe the slow material wave, while the fast variable accounts for the fast acoustic waves. They apply their technique to isentropic and full

Euler and Navier-Stokes, as well as to the isentropic Navier-Stokes-Poisson system.

A slightly different approach is adopted in [37], where the author propose methods based on the flux splitting: the flux is split in two terms, one of which is treated explicitly and the other implicitly.

In a recent paper [2] semi-implicit schemes are constructed. 1D compressible Euler equations, in which acoustic waves are treated implicitly by central central discretization, while the material waves are treated by upwind scheme.

In the present paper we adopt a different strategy. We still treat acoustic waves implicitly and material waves explicitly, however space discretization is obtained by central scheme on staggered grid, thus avoiding the difficulty related to the choice of the exact or approximate Riemann solver to be used when computing the contribution to the numerical flux of the terms which are treated explicitly. This results in a very simple scheme that can be adopted for Euler equations in one and two dimensions. Staggered discretization in space naturally provides second order accuracy, allows a compact stencil in the discretization of the equation for pressure (or energy, according to which variable is chosen as primary unknown). Second order extension in time can be achieved by using globally stiffly accurate IMEX Runge-Kutta (R-K), [9, 8, 7].

The plan of the paper is the following. After the introduction, we define the low Mach number scaling adopted in the paper, and its implications in the isotropic gas dynamics. The next sections are devoted to the construction of the schemes for isentropic gas dynamics in one and two space dimensions: Section 3 deals with first order scheme, while section 4 deals with the extension to second order in time, obtained using globally stiffly accurate IMEX schemes. In Section 5 several numerical tests are presented, and compared with results for the literature. In Section 6 we describe how to extend the schemes to full Euler equations, and in the next Section 7 we present numerical tests on a wide range of problems. In the last section we draw conclusions.

## 2 Low Mach number scaling

We consider the compressible Euler equation for an ideal gas:

$$\begin{cases} \rho_t + \nabla \cdot (\rho \mathbf{u}) = 0 \\ (\rho \mathbf{u})_t + \nabla \cdot (\rho \mathbf{u} \otimes \mathbf{u}) + \nabla p = 0 \\ E_t + \nabla \cdot [(E + p)\mathbf{u}] = 0, \end{cases} \quad (2.1)$$

where $\rho$ is the mass density, $\mathbf{u}$ the velocity of the fluid, $E$ the total energy density per unit volume and $p$ the pressure. System (2.1) is closed by equation of state (EOS), $e = e(\rho, p)$ with $e$ the internal energy, related to total energy density by

$$E = \frac{1}{2}\rho u^2 + \rho e.$$

For polytropic gas we get:

$$e = \frac{1}{\gamma - 1}\frac{p}{\rho}$$

with $\gamma = C_p/C_v > 1$ being the ratio of specific heats of the gas.

In order to describe the low Mach number limit, we rescale the equations considering: $\rho_0, u_0, p_0, x_0, t_0$, with $u_0 = x_0/t_0$ where the dimensionless variables are then given by typical values of the field variables

$$\hat{\rho} = \rho/\rho_0, \quad \hat{u} = u/u_0, \quad \hat{p} = p/p_0, \quad \hat{E} = E/p_0, \quad \hat{x} = x/x_0, \quad \hat{t} = t/t_0.$$

inserting these expressions into the equations (2.1) (and omitting the hat) one obtains the rescaled (non-dimensionalised) compressible Euler equations:

$$\begin{cases} \rho_t + \nabla \cdot (\rho \mathbf{u}) = 0 \\ (\rho \mathbf{u})_t + \nabla \cdot (\rho \mathbf{u} \otimes \mathbf{u}) + \frac{1}{\varepsilon^2} \nabla p = 0 \\ E_t + \nabla \cdot [(E+p)\mathbf{u}] = 0, \end{cases} \quad (2.2)$$

with the equation of state (for a polytropic gas)

$$E = \frac{p}{\gamma - 1} + \frac{\varepsilon^2}{2} \rho |\mathbf{u}|^2$$

where the square of reference Mach number is $\varepsilon^2 = \rho_0 u_0^2/p_0$, strictly speaking, the reference Mach number is $M = \dfrac{u_0}{c_{s_0}} = u_0 \sqrt{\dfrac{\rho_0}{\gamma p_0}} = \dfrac{\varepsilon}{\sqrt{\gamma}}$. This parameter $\varepsilon$ represents a global Mach number characterizing the non dimensionalization but not the local Mach number. System (2.2) is *hyperbolic* and the eigenvalues in direction $\mathbf{n}$ are: $\lambda_1 = \mathbf{u} \cdot \mathbf{n} - c_s/\varepsilon$, $\lambda_2 = \mathbf{u} \cdot \mathbf{n}$, $\lambda_3 = \mathbf{u} \cdot \mathbf{n} - c_s/\varepsilon$ with $c_s = \sqrt{\dfrac{\gamma p}{\rho}}$.

The aim of this work is to construct and analyze new numerical schemes for unsteady compressible flows when the Mach number $\varepsilon$ spans by orders of magnitude.

Compressible flow equations converge to incompressible equations when the Mach number vanishes. This convergence has been rigorously studied mathematically by Klainerman and Majda [27, 29]. When the Mach number is of order one, modern shock capturing methods are able to capture shocks and other complex structures with high numerical resolutions at a reasonable cost.

On the other hand, when we are near the incompressible regime, and Mach number is very small, flows are slow compared with the speed of sound and in such a situation, pressure waves become very fast compared to material waves. In several cases, acoustic waves possess very small energy and they are unimportant near the incompressible regime, then one is not interested in resolving them.

From a numerical point of view, when the Mach number is very small, standard explicit shock-capturing methods require a CFL time restriction dictated by the sound speed $c_s/\varepsilon$ to integrate the system. This leads to the stiffness in time, ([22], [19], [16]) , where the time discretization is constrained by a stability condition given by

$$\Delta t < \Delta x / \lambda_{\max} \approx \mathcal{O}(\varepsilon \Delta x)$$

for small $\varepsilon$ where $\Delta t$ is the time-step, $\Delta x$ the space step and

$$\lambda_{\max} = \max_{\Omega}(|\mathbf{u}| + c_s/\varepsilon)$$

This restriction results in an increasingly large computational time for smaller and smaller $\varepsilon$. The second drawback is due to the excessive numerical viscosity of standard upwind schemes, that scales as $\varepsilon^{-1}$, leading to highly inaccurate solutions. Thus, it is also crucial how the space derivatives are discretized in order to get stability and consistency for the scheme in the incompressible limit (*asymptotic preserving* property).

Our goal in this paper is to develop a numerical scheme that works in all regimes of Mach number for the solution of system (2.2), including both compressible and incompressible regime.

The idea is to design a second order numerical scheme for compressible Euler, whose stability and accuracy is independent of $\varepsilon$, and which is able to capture shocks and discontinuities in the compressible regime, for large $\varepsilon$ and, at the same time, it is a good incompressible solver in the limit regime of vanishing $\varepsilon$. This means that the scheme has to be asymptotic preserving [25], [26], i.e., a numerical scheme which gives a consistent discretization of the compressible Euler equations and in the limit as $\varepsilon \to 0$, with $\Delta x$ and $\Delta t$ fixed, provides a consistent discretization of the incompressible Euler equations. Of course, if a scheme is AP, a uniform accuracy for all range of the parameter $\varepsilon$ is expected.

A key feature of the scheme is the implicit treatment of acoustic waves, while material waves are treated explicitly. Implicit-explicit Euler provides a first order scheme for the Euler equations which filters out the acoustic waves.

For the space discretization, we adopt central schemes on staggered grid similar to the Nessyahu and Tadmor [36] in one space dimension and Jiang and Tadmor in two space dimensions [24], which provided explicit, second order accurate in space and time, shock capturing schemes.

The generalization of NT and JT schemes to higher order was given by the Central Runge-Kutta methods [39]. Here we adopt this approach, coupled with semi-implicit IMEX R-K [6], in order to obtain high order accuracy in time.

The choice of central schemes appears natural because it simplifies flux computation, avoids the introduction of excessive numerical diffusion of upwind discretization, and provides the natural central discretization for the implicit terms.

## 2.1 Isentropic Euler Equations

For sake of clarity, we start considering the isentropic gas dynamics case and successively we extend the results to the case of the full Euler system.

The isentropic Euler equations in d-dimensions, $x \in \Omega \subset \mathbb{R}^d$, $t \geq 0$, are given by:

$$\begin{cases} \rho_t + \nabla \cdot (\rho \mathbf{u}) = 0 \\ (\rho \mathbf{u})_t + \nabla \cdot (\rho \mathbf{u} \otimes \mathbf{u}) + \nabla p(\rho)/\varepsilon^2 = 0, \end{cases} \quad (2.3)$$

where $\rho$ is the density of the fluid, $\mathbf{u}$ is the velocity of the fluid, and $p$ is the pressure. Here we consider a polytropic gas, for which the equation of states take the form: $p(\rho) = C\rho^\gamma$ where $C(s)$ depends on the entropy (which is assumed to be constant) and $\gamma = C_p/C_v$ is the polytropic constant. Here $\varepsilon$ is the dimensionless reference Mach number. As boundary conditions we set $\mathbf{u} \cdot n = 0$ on $\partial\Omega$, or assume $\Omega$ is $\mathbb{T}^d$, i.e. periodic boundary conditions.

Now we recall the classical formal derivation of the incompressible Euler equations from the isentropic compressible Euler system (2.3). We consider an asymptotic expansion ansatz for the following variables:

$$\begin{aligned} \rho(x,t) &= \rho_0(x,t) + \varepsilon^2 \rho_2(x,t) + \cdots, \\ p(x,t) &= p_0(x,t) + \varepsilon^2 p_2(x,t) + \cdots, \\ \mathbf{u}(x,t) &= \mathbf{u}_0(x,t) + \varepsilon^2 \mathbf{u}_2(x,t) + \cdots, \end{aligned} \quad (2.4)$$

we skip the $\mathcal{O}(\varepsilon)$ term because it does not appear in the system equations (2.3). Inserting (2.4) in (2.3), to $\mathcal{O}(\varepsilon^{-2})$ one gets, in the momentum conservation equation (3.2):

$$\nabla p_0 = 0.$$

Therefore, $p_0(x,t) = p_0(t)$, and by $p = p(\rho)$, we have $\rho_0 = \rho_0(t)$, i.e. to lower order in $\varepsilon$ density and pressure are constant in space.

Next, by taking the $\mathcal{O}(\varepsilon^0)$ terms, we have

$$\partial_t \rho_0 + \nabla \cdot (\rho_0 \mathbf{u}_0) = 0 \quad (2.5)$$

$$\partial_t (\rho_0 \mathbf{u}_0) + \nabla \cdot (\rho_0 \mathbf{u}_0 \otimes \mathbf{u}_0) + \nabla p_2 = 0. \quad (2.6)$$

where $p_2 = \lim_{\varepsilon \to 0} \varepsilon^{-2}(p(\rho) - p_0)$ is the hydrostatic pressure. Now, the incompressibility is forced by using the boundary conditions to solve system (2.3) on the domain $\Omega$ with $\mathbf{u} \cdot n = 0$ on $S = \partial\Omega$ or periodic boundary conditions. Because $\rho_0 = \rho_0(t)$ for (2.6) one has:

$$\nabla \cdot \mathbf{u}_0 = -\frac{1}{\rho_0} \frac{d\rho_0}{dt}.$$

Integrating in $\Omega$ one has:

$$-|\Omega| \frac{1}{\rho_0} \frac{d\rho_0}{dt} = \int_\Omega \nabla \cdot \mathbf{u}_0 d\Omega = \int_{\partial\Omega} \mathbf{u}_0 \cdot n dS = 0, \quad (2.7)$$

because of the boundary conditions, therefore $\rho_0 = Const$ (see [22, 16] for more details).

This means that the density satisfies the expansion:

$$\rho(x,t) = \rho_0 + \varepsilon^2 \rho_2(x,t) + \cdots,$$

where $\rho_0$ is a constant of order 1, and then one obtaines $\nabla \cdot \mathbf{u}_0 = 0$. We finally obtain that, for low Mach number, i.e. $\varepsilon \ll 1$, by considering *well-prepared* initial conditions in the sense:

$$\begin{aligned} \rho(x,0) &= \rho_0 + \varepsilon^2 \rho_2(x) + \cdots, \\ \nabla \cdot \mathbf{u}(x,0) &= \mathbf{u}_0(x) + \mathcal{O}(\varepsilon). \end{aligned} \quad (2.8)$$

where $||\mathbf{u}_0||$ is of order 1, such that $\nabla \cdot \mathbf{u}_0 = 0$, the solution $(\rho, \mathbf{u})$ with $p = p(\rho)$ of the isentropic Euler system (2.3) will be close to the solution of the the incompressible Euler system,

$$\begin{aligned} \rho_0 &= Const, \\ \nabla \cdot \mathbf{u}_0 &= 0, \\ \partial_t \mathbf{u}_0 + (\mathbf{u}_0 \cdot \nabla)\mathbf{u}_0 + \frac{\nabla p_2}{\rho_0} &= 0. \end{aligned} \qquad (2.9)$$

We note that, in the low-Mach number model, $p_2$ is the Lagrange multiplier needed to impose the divergence-free constraint: $\nabla \cdot \mathbf{u} = 0$.

Then, taking the divergence of the last equation in (2.9) and using the incompressibility, one obtains

$$-\Delta p_2 = \nabla \cdot (\rho_0 \mathbf{u}_0 \cdot \nabla \mathbf{u}_0) = \nabla^2 : (\rho_0 \mathbf{u}_0 \otimes \mathbf{u}_0). \qquad (2.10)$$

Furthermore, it is possible to derive the pressure wave equation by (2.3), indeed, if we differentiate with respect the time the density equation and subtract it from the divergence of the momentum equation, we obtain

$$\partial_{tt}\rho - \frac{\Delta p(\rho)}{\varepsilon^2} = \nabla^2 : (\rho \mathbf{u} \otimes \mathbf{u}),$$

and to $\mathcal{O}(\varepsilon^0)$, we get (2.10).

Next we propose a numerical scheme that is applicable for all ranges of the Mach number.

# 3 Numerical Schemes

In this section we design a second numerical scheme for compressible Euler that is able to capture the incompressible Euler limit as $\varepsilon \to 0$, i.e. an asymptotic preserving (AP) scheme. Furthermore such a scheme is a conservative and shock capturing for all $\varepsilon$ and has to satisfy a CFL condition which is independent of $\varepsilon$.

The features of our scheme are the following: for the spatial discretization, we use a second-order, non-oscillatory central scheme on a staggered grid in order to simplify the computation of the numerical flux, similar to the ones adopted in [36, 24, 39]. For time integration we start by describing a first order implicit-explicit Euler scheme, while a semi-implicit approach based on IMEX Runge-Kutta methods ([6, 4, 5, 7]) will be presented in Section 4 for higher order generalization.

## 3.1 Isentropic Gas Dynamics: first order method

**1D Model.** For simplicity, we consider the domain $\Omega = [0,1]$, with periodic boundary conditions. Here $\rho$, $u$ and $m$ denote the density, velocity and momentum in one dimension.

Then (2.3) becomes
$$\rho_t + m_x = 0$$
$$m_t + \left(\frac{m^2}{\rho} + \frac{p}{\varepsilon^2}\right)_x = 0 \tag{3.1}$$

The system is closed by $p = \rho^\gamma$. We shall discretize space in a way similar to the NT central scheme (see [36]), i.e., we make use a staggered grid with a uniform spatial mesh $\Delta x = 1/N$, where $N$ is an positive integer and at even time we have $N$ cells of size $\Delta x$, with cell $j$ centered at $x_j = (j - 1/2)\Delta x$, $j = 1, \cdots, N$. While we discretize time by a first order implicit-explicit Euler (5.1): stiff terms will be evaluated at time $t^{n+1}$, while non-stiff terms will be evaluated at time $t^n$.

Then, integrate the equation on a staggered grid, from time $t^n = n\Delta t$, $n = 0, 1. \cdots$, to $t^{n+1}$ (see Figure 3.1) we obtain the first order semi-implicit scheme:

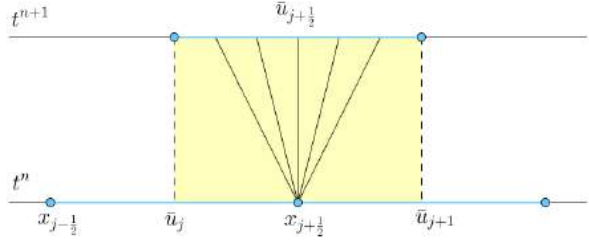

Figure 3.1: Staggered grid form $t^n$ to $t^{n+1}$.

$$\bar{\rho}_{j+1/2}^{n+1} = \bar{\rho}_{j+1/2}^n - \frac{\Delta t}{\Delta x}(m_{j+1}^{n+1} - m_j^{n+1})$$
$$\bar{m}_{j+1/2}^{n+1} = \bar{m}_{j+1/2}^n - \frac{\Delta t}{\Delta x}(f_{j+1}^n - f_j^n) - \frac{\Delta t}{\varepsilon^2 \Delta x}(p_{j+1}^{n+1} - p_j^{n+1}) \tag{3.2}$$

where $f_j^n = (\bar{m}_j^n)^2/\bar{\rho}_j^n$. We note that second order in space is obtained by standard reconstruction adopted in Nessyahu-Tadmor scheme (see for details [36]), e.g.

$$\bar{\rho}_{j+1/2} = \frac{\bar{\rho}_{j+1} + \bar{\rho}_j}{2} + \frac{1}{8}(\hat{D}_x\rho_j - \hat{D}_x\rho_{j+1}),$$

with $\Delta x \hat{D}_x \rho_j$ a first order approximation of the first derivative on cell $j$, for example,

$$\hat{D}_x\rho_j = \text{MM}(\rho_j - \rho_{j-1}, \rho_{j+1} - \rho_j)$$

or

$$\hat{D}_x\rho_j = \text{MM}\left(\theta(\rho_j - \rho_{j-1}), \frac{\rho_{j+1} - \rho_{j-1}}{2}, \theta(\rho_{j+1} - \rho_j)\right) \tag{3.3}$$

where MM denotes the minimod function and $\theta \in [1, 2]$.

The flux term appearing in the first equation of (3.2), computed at cell center, are defined as:
$$m_j^{n+1} = m_j^* - \frac{\Delta t}{\varepsilon^2 \Delta x} D_x p_j^{n+1}$$
where $m_j^* = m_j^n - \frac{\Delta t}{\Delta x}\hat{D}_x f_j^n$, and $D_x$ denotes central difference on the cell: $D_x p_j^{n+1} = (p_{j+\frac{1}{2}}^{n+1} - p_{j-\frac{1}{2}}^{n+1})$. Using such equation, and substituting it into the density equation for $\bar{\rho}_{j+1/2}^{n+1}$ one gets an equation of the form:

$$\bar{\rho}_{j+1/2}^{n+1} - \frac{\Delta t^2}{\varepsilon^2 \Delta x^2} D_x^2 p_{j+1/2}^{n+1} = \rho_{j+1/2}^* \tag{3.4}$$

where $D_x^2 h_j = h_{j+1} - 2h_j + h_{j-1}, \forall h_j$, is the usual three point discrete Laplacian and

$$\rho_{j+1/2}^* = \bar{\rho}_{j+1/2}^n - \frac{\Delta t}{\Delta x}(m_{j+1}^* - m_j^*)$$

denotes quantities that can be computed explicitly (in a conservative way). When we use the (second order) approximation $p_{j+1/2}^{n+1} = p(\bar{\rho}_{j+1/2}^{n+1})$, (3.4) becomes a non linear equation for the new density on the staggered mesh.

One possible way to simplify the solution if the system is to use $p$ as unknown and considering $\rho = \rho(p)$, then we get

$$\left(p_{j+1/2}^{n+1}\right)^{1/\gamma} = \rho_{j+1/2}^* - \frac{\Delta t^2}{\varepsilon^2 \Delta x^2} D_x^2 p_{j+1/2}^{n+1}. \tag{3.5}$$

In this case the nonlinearity is in the diagonal of the system, and the linear equation for each time step for the unknown $p_{j+1/2}^{n+1}$ can be solved by few iterations of Newton's method.

On the other hand, if we approximate the Laplace operator $(p^{n+1})_{xx}$ discretized in (3.4) by the approximation

$$(p(\rho^{n+1}))_{xx} \approx (p(\rho^n)' \rho_x^{n+1})_x$$

a semi-implicit approach is involved and this implies to solve a linear system in (3.4).

**Remark.** We observe that the use of a staggered grid allows a compact discrete Laplacian in the implicit equation (3.4) for the new density ((3.5) for the pressure). This property is important, ant it is one of the main reasons, many solvers for incompressible Euler and Navier-Stokes equations in primitive variables are discretized on a staggered grid (see for example [42], Sec. 6.3).

**2D Model.** Now, we consider the domain $\Omega = [0,1] \times [0,1]$ with periodic boundary conditions. The Euler equations for isentropic case in 2D of the gas dynamics are given

by
$$\rho_t + \nabla \cdot \mathbf{m} = 0$$

$$\frac{\partial \overset{1}{m}}{\partial t} + \nabla \cdot \left(\frac{\overset{1}{m}\mathbf{m}}{\rho}\right) + \frac{\partial_x p}{\varepsilon^2} = 0$$
$$\frac{\partial \overset{2}{m}}{\partial t} + \nabla \cdot \left(\frac{\overset{2}{m}\mathbf{m}}{\rho}\right) + \frac{\partial_y p}{\varepsilon^2} = 0 \quad (3.6)$$

where $\mathbf{m} = (\overset{1}{m}, \overset{2}{m})^T$. The system is closed by $p = C\rho^\gamma$. We choose $C = 1$.

First of all we discretize (3.6) time by an explicit-implicit Euler scheme:

$$\rho^{n+1} = \rho^n - \Delta t \nabla \cdot \mathbf{m}^{n+1},$$
$$\overset{1}{m}{}^{n+1} = \overset{1}{m}{}^n - \Delta t \nabla \cdot \left(\frac{\overset{1}{m}{}^n \mathbf{m}^n}{\rho^n}\right) - \Delta t \frac{\partial_x p^{n+1}}{\varepsilon^2}, \quad (3.7)$$
$$\overset{2}{m}{}^{n+1} = \overset{2}{m}{}^n - \Delta t \nabla \cdot \left(\frac{\overset{2}{m}{}^n \mathbf{m}^n}{\rho^n}\right) - \Delta t \frac{\partial_y p^{n+1}}{\varepsilon^2}.$$

Let

$$\overset{k}{m}{}^* = \overset{k}{m}{}^n - \Delta t \nabla \cdot \left(\frac{\overset{k}{m}{}^n \mathbf{m}^n}{\rho^n}\right), \quad k = 1, 2, \quad (3.8)$$

be the explicit part of the second and third equation in (3.7). Then (3.7) becomes

$$\rho^{n+1} = \rho^n - \Delta t \nabla \cdot \mathbf{m}^{n+1},$$
$$\overset{1}{m}{}^{n+1} = \overset{1}{m}{}^* - \Delta t \frac{\partial_x p^{n+1}}{\varepsilon^2}, \quad (3.9)$$
$$\overset{2}{m}{}^{n+1} = \overset{2}{m}{}^* - \Delta t \frac{\partial_y p^{n+1}}{\varepsilon^2},$$

inserting the expression $\mathbf{m}^{n+1}$ into the first equation for $\rho$ in we obtain

$$\rho^{n+1} = \rho^* + \frac{\Delta t^2}{\varepsilon^2} \nabla^2 p^{n+1}, \quad (3.10)$$

where

$$\rho^* = \rho^n - \Delta t \nabla \cdot \mathbf{m}^*. \quad (3.11)$$

Now we discretize space in a way similar to the JT central scheme (see [24]), i.e., we make use a staggered grid with a uniform spatial mesh with $\Delta x = 1/N$, $\Delta y = 1/N$, with $N$ a positive integer, the grid points defined as $x_i = i\Delta x$, $i = 0, 1, \cdots, N$, and $y_j = j\Delta y$, $j = 0, 1, \cdots, N$. In the JT central scheme based on staggered grid, system (3.6) is integrated on the staggered control volume: $C_{i+\frac{1}{2}, j+\frac{1}{2}} \times [t^n, t^{n+1})$ with $C_{i+\frac{1}{2}, j+\frac{1}{2}} := I_{i+\frac{1}{2}} \times J_{j+\frac{1}{2}}$ centered around $(x_{i+\frac{1}{2}}, y_{j+\frac{1}{2}})$, with $I_{i+\frac{1}{2}} = [x_i, x_{i+1}]$ and $J_{j+\frac{1}{2}} = [y_j, y_{j+1}]$ (see fig. 3.1).

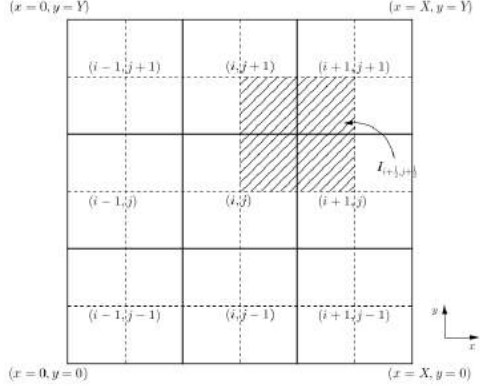

Figure 3.2: Control Volume: $C_{i+1/2,j+1/2}$

We assume we know cell averages $\bar{\rho}^n_{ij}$, $\bar{\mathbf{m}}^n_{ij}$ on the initial grid at even $n$, and we want to compute the field variables $\bar{\rho}^{n+1}_{i+1/2,j+1/2}$, $\bar{\mathbf{m}}^{n+1}_{i+1/2,j+1/2}$ on the staggered grid.

The scheme works as follows.

1. Compute $\mathbf{m}^*$ on the staggered grid discretizing Eq. (3.8):

$$\overset{k}{m}{}^*_{i+1/2,j+1/2} = \overset{k}{\bar{m}}{}^n_{i+1/2,j+1/2} - \Delta t \mathcal{D}\left(\frac{\overset{k}{\bar{m}}{}^n \bar{\mathbf{m}}^n}{\bar{\rho}^n}\right)_{i+1/2,j+1/2}, \quad k = 1, 2. \qquad (3.12)$$

2. Evaluating $\bar{\rho}^{n+1}_{i+\frac{1}{2},j+\frac{1}{2}}$ on the staggered grid discretizing the first Eq. in (3.6)

$$\bar{\rho}^{n+1}_{i+\frac{1}{2},j+\frac{1}{2}} = \bar{\rho}^n_{i+1/2,j+1/2} - \Delta t \mathcal{D}\left(\bar{m}^{n+1}_{i+1/2,j+1/2}\right), \qquad (3.13)$$

compute $\rho^*$ (3.11) by:

$$\rho^*_{i+1/2,j+1/2} = \bar{\rho}^n_{i+1/2,j+1/2} - \Delta t \mathcal{D}\left(\bar{m}^*_{i+1/2,j+1/2}\right). \qquad (3.14)$$

where the staggered cell averages density is computed as [24],

$$\begin{aligned}
\bar{\rho}_{i+\frac{1}{2},j+\frac{1}{2}} &= \tfrac{1}{4}(\bar{\rho}_{i,j} + \bar{\rho}_{i+1,j} + \bar{\rho}_{i,j+1} + \bar{\rho}_{i+1,j+1}) \\
&+ \tfrac{1}{16}(\rho'_{i,j} - \rho'_{i+1,j} + \rho'_{i,j+1} - \rho'_{i+1,j+1})\Delta x \\
&+ \tfrac{1}{16}(\rho^{\backslash}_{i,j} - \rho^{\backslash}_{i,j+1} + \rho^{\backslash}_{i+1,j} - \rho^{\backslash}_{i+1,j+1})\Delta y
\end{aligned}$$

whith $\rho'_{i,j}$ and $\rho^{\backslash}_{i,j}$ a first order approximation of the first derivative on cell $(i,j)$ (in this paper we use minmod in most cases).

3. Solve the main linear elliptic equation (3.10):

$$\bar{\rho}^{n+1}_{i+\frac{1}{2},j+\frac{1}{2}} = \bar{\rho}^{*}_{i+1/2,j+1/2} + \frac{\Delta t^2}{\varepsilon^2} L p^{n+1}_{i+1/2,j+1/2}. \tag{3.15}$$

4. Compute $\bar{\mathbf{m}}^{n+1}$:

$$\overset{k}{m}{}^{n+1}_{i+1/2,j+1/2} = \overset{k}{\bar{m}}{}^{*}_{i+1/2,j+1/2} - \frac{\Delta t}{\varepsilon^2} D_k p^{n+1}_{i+1/2,j+1/2}, \quad k = 1, 2.$$

where $D_k$ is the classical central difference approximation of the space derivative in the $k$-th direction

$$D_1 h_{ij} = \frac{h_{i+1,j} - h_{i-1,j}}{2\Delta x}, \quad D_2 h_{ij} = \frac{h_{i,j+1} - h_{i,j-1}}{2\Delta y}.$$

.

The operators $\mathcal{D}$ and $L$ are defined as:

$$\mathcal{D}\mathbf{m}_{i+1/2,j+1/2} = \frac{\overset{1}{m}_{i+1,j} - \overset{1}{m}_{i,j}}{2\Delta x} + \frac{\overset{1}{m}_{i+1,j+1} - \overset{1}{m}_{i,j+1}}{2\Delta x} + \frac{\overset{2}{m}_{i,j+1} - \overset{2}{m}_{i,j}}{2\Delta y} + \frac{\overset{2}{m}_{i+1,j+1} - \overset{2}{m}_{i+1,j}}{2\Delta y},$$

and

$$L p_{ij} = (p_{i+1,j} - 2p_{i,j} + p_{i-1,j})/\Delta x^2 + (p_{i,j+1} - 2p_{i,j} + p_{i,j-1})/\Delta y^2.$$

$L p_{ij}$ is the discrete Laplacian in 2D.

Note that when we use the relation $p^{n+1}_{i+1/2,j+1/2} = p(\bar{\rho}^{n+1}_{i+1/2,j+1/2})$, (3.15) becomes a non linear equation for the new density in the staggered mesh, but if one takes $p$ as unknown and considering $\rho = \rho(p)$ then we solve a non linear system by a Newton's method where the nonlinearity is in the diagonal of the system as in the 1D case. However, linearizing the operator $\nabla^2 p^{n+1}$ as $\nabla \cdot (p'(\rho^n)\nabla \rho^{n+1})$, this implies again to solve a linear system in (3.15) for the new density on the staggered mesh.

**Remark.** Notice that for even $n$ the discrete operators of quantities at time $t^n$ act on discrete field values with integer indices, while at time $t^{n+1}$, they act on field values with half-odd integers, thus maintaining the operatorr compact. A similar procedure is adopted when going from $n$ add to $n+1$, thus obtaining discrete field variables again on the original mesh.

## 3.2 Asymptotic preserving (AP) property

In order to show the AP property of our scheme we should demonstrate that such a scheme in the low Mach number limit, i.e. as $\varepsilon \to 0$, provides a consistent discretization of the incompressible Euler equation (2.9) with spatial and temporal steps fixed.

For this analysis, by (2.8), we consider the following *well-prepared* initial conditions:

$$\begin{aligned} \rho^n_{ij} &= \rho^n_{0,ij} + \varepsilon^2 \rho^n_{2,ij} + \cdots, \\ m^n_{ij} &= (\rho u)^n_{ij} = (\rho u)^n_{0,ij} + \varepsilon^2 (\rho u)^n_{2,ij} + \cdots. \end{aligned} \tag{3.16}$$

By substituting this ansatz (3.16) into the scheme (3.12)-(3.14)-(3.15), we get for $\mathcal{O}(\varepsilon^{-2})$ terms
$$Lp^{n+1}_{0,i+1/2,j+1,2} = 0,$$
and by the fact that $p(\rho) = \rho^\gamma$, one gets space independence for the leading order density, i.e., $\rho^{n+1}_{0,ij} = \rho^{n+1}_{(0)}$ for all $i,j$, then it is constant in space but not necessarily in time. Now, the $\mathcal{O}(1)$ equation for the density, from (3.13), is given by:

$$\rho^{n+1}_{0,i+1/2,j+1/2} = \bar{\rho}^n_{0,i+1/2,j+1/2} - \Delta t \mathcal{D}\left(\bar{m}^{n+1}_{0,i+1/2,j+1/2}\right), \tag{3.17}$$

Now considering the analog discretization of the integral (2.7), and periodic boundary conditions, we summing the previous equation over all $i,j$ and by the fact that $\rho^{n+1}_{0,ij} = \rho^{n+1}_0$ for all $i,j$, we get

$$\bar{\rho}^{n+1}_{(0)} = \frac{1}{N}\sum_{ij}\bar{\rho}^n_{(0),i+1/2,j+1/2} := \bar{\rho}_{(0)}. \tag{3.18}$$

where $N$ is the total number of grid points. This previous relation says that the new density is a constant and equal to the average value of the density at the previous time step. Furthermore, by the consistency with the initial date (2.8), the density, in the previous step $rho^n_0$ is constant, then the two quantities in (3.18) coincide, so that the density is also constant in time, just as in the continuous case in the incompressible regime.

Now, using this result in Eq (3.17), the density terms cancel out and we get, $\mathcal{D}\left(\overset{k}{m}{}^{n+1}_{(0),i+1/2,j+1/2}\right) = 0$ with $\overset{k}{m} = \rho\overset{k}{u}$ for $k = 1, 2$. Then, we obtain the discrete incompressibility condition for the vector velocity, i.e.,

$$\mathcal{D}\mathbf{u}^{n+1}_{(0),i+1/2,j+1/2} = 0. \tag{3.19}$$

Furthermore, we observe that assuming that the initial velocity field is incompressible in the (2.8), i.e. $\nabla \cdots \mathbf{u}_0 = 0$, this means that at the previous time step $\nabla \cdots \mathbf{u}^n_0 = 0$ and then (3.19) is satisfied for all $n$. Now, by this fact follows that $\mathcal{D}\mathbf{m}^n_{(0)i+1/2,j+1/2} = 0$. Now we derive an equation for the pressure term $p_{2,ij}$. Then, considering the $\mathcal{O}(1)$ terms in the linear elliptic equation (3.15):

$$\bar{\rho}^{n+1}_{(0),i+1/2,j+1/2} = \rho^*_{(0),i+1/2,j+1/2} + \Delta t^2 Lp^{n+1}_{(2),i+1,2,j}, \tag{3.20}$$

we get, by inserting (3.12) in (3.14) and some algebraic computation:

$$-Lp^{n+1}_{(2),i+1,2,j+1/2} = \mathcal{D}\left(\mathcal{D}(\rho_{(0)}\mathbf{u}_{(0)} \otimes \mathbf{u}_{(0)})_{i+1/2,j+1/2}\right)_{i+1/2,j+1/2}. \tag{3.21}$$

Note that (3.21) is the discretization of the equation (2.10).

Finally, by the knowledge of $\bar{\rho}_{(0)}$ is independent of the time, the $\mathcal{O}(1)$ term of the momentum equation becomes:

$$\bar{\rho}_{(0)} \left( \frac{\bar{\mathbf{u}}^{n+1}_{(0),i+\frac{1}{2},j+\frac{1}{2}} - \bar{\mathbf{u}}^{n}_{(0),i+\frac{1}{2},j+\frac{1}{2}}}{\Delta t} \right) = \mathcal{D}\left(\bar{\rho}_{(0)}\mathbf{u}_{(0)} \otimes \mathbf{u}_{(0)}\right)_{i+1/2,j+1/2} - \mathcal{D}p^{n+1}_{(2),i+1/2,j+1/2}. \tag{3.22}$$

Thus, (3.18)-(3.19) and (3.22) represent a discretization of the incompressible Euler equation (2.9). Therefore, the two-dimensional scheme is AP.

## 4 Extension to a second order

In this section we present a second order scheme for the isentropic Euler equations.

The approach will be a combination of central Runge-Kutta methods for conservation laws (CRK, [39]) and semi-implicit IMEX schemes developed in [6]. In some recent papers [7, 6, 4, 5] a very effective semi-implicit technique has been introduced for the numerical solution of systems nonlinear hyperbolic systems. The method is based on explicit-implicit Runge-Kutta methods (IMEX R-K) and the proposed schemes are stable, linearly implicit, and can be designed up to any order of accuracy.

IMEX RK schemes can be represented by a pair Butcher arrays given by:

$$\begin{array}{c|c} \tilde{c} & \tilde{A} \\ \hline & \tilde{b}^T \end{array}, \quad \begin{array}{c|c} c & A \\ \hline & b^T \end{array}.$$

where the $s \times s$ low triangular matrices $\tilde{A} = (\tilde{a}_{ij})$ ($\tilde{a}_{ij} = 0$ for all $j \geq i$), and $A = (a_{ij})$ ($a_{ij} = 0$ for all $j > i$) are the matrices of the explicit and implicit parts of the scheme, respectively, while the vectors $\tilde{b} = (\tilde{b}_1, \cdots, \tilde{b}_s)$, $b = (b_1, \cdots, b_s)$, $\tilde{c} = (\tilde{c}_1, \cdots, \tilde{c}_s$, and $c = (c_1, \cdots, c_s)$, are $s$-dimensional vectors or real coefficients, which $\tilde{c}$ and $c$ given by the usual relations

$$\tilde{c}_i = \sum_{i=1}^{i-1} \tilde{a}_{ij}, \quad c_i = \sum_{i=1}^{i} a_{ij}, \quad i = 1, \cdots, s.$$

Furthermore, in order to reduce the computational cost for the implicit part of the IMEX R-K scheme, we consider a particular class of implicit R-K methods called *diagonally implicit R-K*, [23].

Before applying the technique developed in [6], we write system (3.6) in the form of an autonomous system

$$\frac{d\mathbf{U}}{dt} = \mathcal{H}(\mathbf{U}, \mathbf{U}). \tag{4.1}$$

with $\mathcal{H} : \mathbb{R}^m \times \mathbb{R}^m \to \mathbb{R}^m$ is a sufficiently regular mapping.

We assume that the dependence on the first argument of $\mathcal{H}$ is non-stiff, while the second argument is stiff, [6]. We further emphasize such dependence by using an asterisk to denote explicit variables: $\mathcal{H} = \mathcal{H}(\mathbf{U}^*, \mathbf{U})$.

In one space dimension, one has $\mathbf{U} = (\rho, m)^T$, and $\mathcal{H}(\mathbf{U}^*, \mathbf{U})$ is:

$$\mathcal{H}(\mathbf{U}^*, \mathbf{U}) = \begin{pmatrix} -m_x \\ -\left(\dfrac{m^2}{\rho}\right)^*_x - \dfrac{p(\rho)_x}{\varepsilon^2} \end{pmatrix}. \tag{4.2}$$

while in 2D one has $\mathbf{U} = (\rho, \mathbf{m})^T$, wtih $\mathbf{m} = (\overset{1}{m}, \overset{2}{m})$ and in this case the function $\mathcal{H}(\mathbf{U}^*, \mathbf{U})$ is:

$$\mathcal{H}(\mathbf{U}^*, \mathbf{U}) = \begin{pmatrix} -\nabla \cdot \mathbf{m} \\ -\nabla \cdot \left(\dfrac{\overset{1}{m}\mathbf{m}}{\rho}\right)^* - \dfrac{\partial_x p(\rho)}{\varepsilon^2} \\ -\nabla \cdot \left(\dfrac{\overset{2}{m}\mathbf{m}}{\rho}\right)^* - \dfrac{\partial_y p(\rho)}{\varepsilon^2} \end{pmatrix}. \tag{4.3}$$

We limit the description of the second order extension to the one dimensional case. The same discretization in time can be applied to the 2D case in a similar way.

Explicit CRK schemes on a staggered grid for a system of conservation laws $u_t + f(u)_x = 0$, can be described as follows, [39].

- **Predictor**. Given cell averages at time $t^n$: $\{\bar{u}_j^n\}$, one computes the stage values as

$$u_j^{(k)} = u_j^n - \frac{\Delta t}{\Delta x} \sum_{\ell=1}^{k-1} a_{k\ell} D_x f(u_j^{(\ell)}), \quad k = 1, \cdots, s,$$

where $\Delta x D_x$ denotes a consistent discrete derivative in space, and $a_{k\ell}$ are the RK coefficients.

- **Corrector**. The numerical solution $\bar{u}_{j+1/2}$ is then recovered as

$$\bar{u}_{j+1/2}^{n+1} = \bar{u}_{j+1/2}^n - \frac{\Delta t}{\Delta x} \sum_{k=1}^{s} b_k \left( f(u_{j+1}^{(k)}) - f(u_{j+1}^{(k)}) \right).$$

where $\bar{u}_{j+1/2}^n$ is obtained integrating the reconstruction of $u^n(x)$ in the staggered cell $I_{j+1/2}$. One might think that the approach to problems containing very stiff terms could be generalized by computing the stage values using IMEX schemes, with an $L$-stable implicit part. However, such *naive* generalization does not work in practice since it gives a CFL type stability restriction comparable with the the one of explicit CRK schemes. We show this by applying the simplest first order scheme with implicit Euler predictor to the linear convection equation $u_t + u_x = 0$.

In this case the scheme (which is first order in space and time) becomes:

- **Predictor**: $u_j^{n+1} = \bar{u}_j - \dfrac{\Delta t}{2\Delta x}(u_{j+1}^{n+1} - u_{j-1}^{n+1})$

- **Corrector**: $\bar{u}_{j+1/2}^{n+1} = \dfrac{\bar{u}_j + \bar{u}_{j+1}}{2} - \dfrac{\Delta t}{\Delta x}(u_{j+1}^{n+1} - u_j^{n+1})$

Looking for solutions of the form: $\bar{u}_j^n = \rho^n e^{ijkh}$, where $k$ denotes the Fourier node and $h = \Delta x$, one obtains the following expression for the amplification factor:

$$\rho = \cos(\xi/2) - \frac{2ic\sin(\xi/2)}{1 + ic\sin\xi},$$

where $\xi = kh$, $c = \Delta t/\Delta x$. Then we get

$$|\rho|^2 = \frac{\cos^2(\xi/2) + 4c^2\sin^6(\xi/2)}{1 + c^2\sin^2\xi} = \frac{\mathcal{N}}{\mathcal{D}}.$$

Let $\mathcal{F} \equiv \mathcal{D} - \mathcal{N} = \sin^2(\xi/2)\left(1 - 4c^2(\sin^4(\xi/2) - \cos^2(xi/2))\right)$. Then $|\rho| \leq 1$ iff $\mathcal{F} \geq 0$. This condition is guaranteed $\forall \xi \in \mathbb{R}$ iff $c \leq 1/2$.

Stability analysis can be performed for other $L$-stable scheme, with similar outcomes. In view of the above results we adopt this approach for the computation of the first $s-1$ stages, at the unstaggered cell centers, however, the last stage $s$ has to be implicit *at the level of the numerical solution.* This can be achieved by imposing that the numerical solution is automatically obtained as the last stage of the scheme directly computed at time $t^{n+1}$ on the staggered cell. In the simple example above, this means

$$\bar{u}_{j+1/2}^{n+1} = \frac{\bar{u}_j + \bar{u}_{j+1}}{2} - \frac{\Delta t}{\Delta x}(\bar{u}_{j+1}^{n+1} - \bar{u}_j^{n+1}) = \frac{\bar{u}_j + \bar{u}_{j+1}}{2} - \frac{\Delta t}{2\Delta x}(u_{j+3/2}^{n+1} - u_{j-1/2}^{n+1}),$$

which is unconditionally stable, i.e., he amplification factor is given by

$$|\rho|^2 = \cos(\xi/2)(1 + c^2\sin^2\xi)^{-1} < 1, \forall \xi.$$

Implicit schemes in which the numerical solution coincides with the last stage are called *stiffly accurate* (SA) (see [23]). IMEX R-K schemes with the same property are called *globally stiffly accurate* (GSA):

**Definition:** We say that an IMEX R-K scheme is *globally stiffly accurate* if it is SA, i.e. $b^T = e_s^T A$, and $\tilde{b}^T = e_s^T \tilde{A}$, with $e_s = (0, \ldots, 0, 1)^T$, and $c_s = \tilde{c}_s = 1$. (see [9, 7, 8]).

In view of the above consideration, we shall obtain high order accuracy in time by adopting GSA IMEX R-K schemes where the first $s-1$ stages are obtained at unstaggered cell center by a predictor (not necessarily conservative), while the numerical solution is obtained by a conservative corrector on the staggered cell.

We now apply the semi-implicit IMEX R-K schemes developed in [6] to system (4.1) with (4.2). The algorithm for the semi-implicit IMEX R-K schemes on staggered cells, (SI-IMEX-RK-STAG) can be conveniently written with two-step, as follows:

1. **Prediction step.** Compute the internal stages:
$\mathbf{U}_j^{*(1)} = \mathbf{U}_i^n$
for $k = 1, \cdots, s-1$

$$\mathbf{U}_j^{*(k)} = \mathbf{U}_j^n - \frac{\Delta t}{\Delta x} \sum_{\ell=1}^{k-1} \tilde{a}_{k,\ell} K_j^{(\ell)}, \tag{4.4a}$$

$$\mathbf{U}_j^{(k)} = \mathbf{U}_j^n - \Delta t \sum_{\ell=1}^{k} a_{k+1,\ell} K_j^{(\ell)}, \tag{4.4b}$$

with the SI-IMEX R-K fluxes:

$$K_j^{(\ell)} = \begin{pmatrix} \hat{D}_x m_j^{(\ell)} \\ \hat{D}_x \left(\frac{m_j^2}{\rho_j}\right)_j^{*(\ell)} + \frac{1}{\varepsilon^2} D_x p(\rho_j^{(\ell)}) \end{pmatrix} \tag{4.5}$$

2. **correction step.** Update the cell averages on the staggered grid:

$$\bar{\mathbf{U}}_{i+1/2}^{n+1} = \bar{\mathbf{U}}_{i+1/2}^n - \frac{\Delta t}{\Delta x} \sum_{\ell=1}^{s} a_{s,\ell} \Delta F_{j+1/2}^{(\ell)}, \tag{4.6}$$

where

$$\Delta F_{j+\frac{1}{2}}^{(\ell)} \equiv F(U_{j+1}^{*(\ell)}, U_{j+1}^{(\ell)}) - F(U_j^{*(\ell)}, U_j^{(\ell)}), \quad \ell = 1, \cdots, s-1$$

and

$$\Delta F_{j+1/2}^{(s)} \equiv F(U_{j+1}^{*(s)}, U_{j+1}^{(s)}) - F(U_j^{*(s)}, U_j^{(s)}).$$

Here

$$F(U^*, U) = \left(m, \left(\frac{m^2}{\rho}\right)^* + \frac{p(\rho)}{\varepsilon^2}\right)^\top$$

.

**Remark.** Throughout the paper we limit to second order accuracy in space and time. In this case one can use $\mathbf{U}_i^n = \bar{\mathbf{U}}_i^n$ in (4.4). The techniques described above can be adopted to construct high order schemes. This requires: i) high order GSA IMEX ([9, 8, 6]), ii) point-wise reconstruction of $\mathbf{U}_i^n$ from $\bar{\mathbf{U}}_i^n$, iii) high order reconstruction (such as WENO [40]) for the terms which are treated explicitly, iv) a high order approximation of the operator appearing in the implicit terms. This research is being currently carried out in [3].

# 5   Numerical Results for Isentropic gas dynamics.

In this section we present the performances of the proposed first and second IMEX-RK-STAG scheme. We test our new schemes presenting several numerical test cases in one and two space dimensions and show that the schemes are accurate for a wide range of values of the Mach number. The schemes run for values of the parameter $\varepsilon$ ranging from compressible to incompressible flows. For all the numerical tests we give well prepared initial values and adopt periodic boundary conditions. Finally, several convergence tests permit to observe the correct second order accuracy of our scheme both in compressible and incompressible regimes.

In all our tests, we used a second order reconstruction (3.3) with $\theta = 1.5$, and we consider the following GSA IMEX R-K schemes.

- First order GSA Euler IMEX scheme [1]:

$$
\begin{array}{c|cc} 0 & 0 & 0 \\ 1 & 1 & 0 \\ \hline & 1 & 0 \end{array} \qquad \begin{array}{c|cc} 0 & 0 & 0 \\ 1 & 0 & 1 \\ \hline & 0 & 1 \end{array} \tag{5.1}
$$

- Second order GSA IMEX R-K scheme:

$$
\begin{array}{c|ccc} 0 & 0 & 0 & 0 \\ c & c & 0 & 0 \\ 1 & 1-1/(2c) & 1/(2c) & 0 \\ \hline & 1-1/(2c) & 1/(2c) & 0 \end{array} \qquad \begin{array}{c|ccc} 0 & 0 & 0 & 0 \\ c & 0 & c & 0 \\ 1 & 0 & 1-\gamma & \gamma \\ \hline & 0 & 1-\gamma & \gamma, \end{array} \tag{5.2}
$$

where $\gamma = (c-1/2)/(c-1)$. In this family of schemes, the choice of $c > 1$ guarantees that $\gamma > 0$ and the weights of the explicit part are non negative. Furthermore, the scheme has an $L$-stable implicit part and an $I$-stable explicit part, [23]. Note that the choice $c = 2.25$ guarantees a small constant of the error in the implicit scheme, so this is the value that we adopt in our scheme. Finally, note that choosing $c = 1 - 1/\sqrt{2}$ we get the classical DIRK(2,2,2) IMEX scheme, [1].

## 5.1 Example 1: Riemann problem

We consider the following initial data, [19]

$$
\begin{cases} \rho(x,0) = 1.0, & m(x,0) = 1 - \frac{\varepsilon^2}{2}, & x \in [0, 0.2] \bigcup [0.8, 1], \\ \rho(x,0) = 1 + \varepsilon^2, & m(x,0) = 1, & x \in (0.2, 0.3], \\ \rho(x,0) = 1, & m(x,0) = 1 + \frac{\varepsilon^2}{2}, & x \in (0.3, 0.7], \\ \rho(x,0) = 1 - \varepsilon^2, & m(x,0) = 1, & x \in (0.7, 0.8], \end{cases} \tag{5.3}
$$

This example consists of several Riemann problems. We choose $p(\rho) = \rho^2$, final time $T = 0.05$ and periodic boundary conditions. In Figure 5.1 we report the solutions for the density on the left and for the momentum on the right. In this example we choose $\varepsilon = 0.8, 0.3, 0.05$. We note that for large $\varepsilon$, shocks ands contact discontinuities appear. We performed also a calculation with $\varepsilon^{-4}$ and the numerical solution quickly relaxes to the trivial solution $\rho = 1, m = 1$, so we omitted to report the figure. In this test we choose the time step $\Delta t$ as follows:

$$\Delta t = CFL_{Imp} \frac{\Delta x}{max|\mathbf{u}_j| + \tilde{c}_j} \qquad (5.4)$$

where $\tilde{c}_j = \sqrt{\gamma p_j/\rho_j} \cdot \min(1, 1/\varepsilon)$, $\forall \varepsilon$ and it is reported in the figures. Note that for $\epsilon > 1$ this is equivalent to the classical CFL condition, but for smaller values of $\epsilon$ the condition is less restrictive.

For each test we monitor the classical Courant number

$$CFL = \frac{\lambda_{\max} \Delta t}{\Delta x}$$

where $\lambda_{\max} = \max_j(|\mathbf{u}_j| + c_j/\varepsilon)$, $c_j = \sqrt{\gamma p_j/\rho_j}$.

Standard explicit schemes on staggered grid require a classical Courant number $CFL \leq 1/2$, while for our semi-implicit scheme $CFL$ can be quite large than 1 for small Mach numbers. Because only the material wave is treated explicitly, we expect the stability restriction for our semi-implicit schemes is given by $CFL_u \leq 1/2$ where

$$CFL_u = \frac{\max |\mathbf{u}_j| \Delta t}{\Delta x}.$$

This is in practice the condition that we adopt for small values of $\epsilon$. We prefer to use condition (5.4) in order to avoid excessively large values of the time step when the velocity of the gas is too small. A formal stability analysis of the scheme is under way.

Finally, a reference solution is computed with $\Delta x = 1/500$, and $\Delta t = 1/20000$.

We observe that for moderate values of $\varepsilon$ the numerical solutions are close to the reference one, while for small values of $\varepsilon$ the acoustic waves are dumped out, still the numerical solution maintains stability.

## 5.2 Example 2: Convergence test.

Here we verify the temporal and spatial order of accuracy of the scheme in compressible and incompressible regimes in one dimension. In order to do that we consider Equations (3.1), take as computational domain $\Omega = [-2.5, 2.5]$ with initial conditions

$$u_0(x) = \sin\left(\frac{2\pi x}{L}\right), \quad \rho_0(x) = \left(1 + \frac{(\gamma-1)u(x)_0}{2c}\right)^{\frac{2}{\gamma-1}}, \quad p_0 = \rho_0^\gamma$$

where $\gamma = 2, c = \sqrt{\gamma}/\varepsilon, L = 5$ and final time $T = 0.3$, $\Delta t$ is given by (5.4) with $CFL_{Imp} = 0.45$. Density errors for different values of the Mach number are listed in Table 5.1. Since we used a staggered grid, we compute the experimental order of convergence (EOC) by the formula

$$EOC := log_2\left(\frac{e_N}{e_{2N}}\right) \qquad (5.5)$$

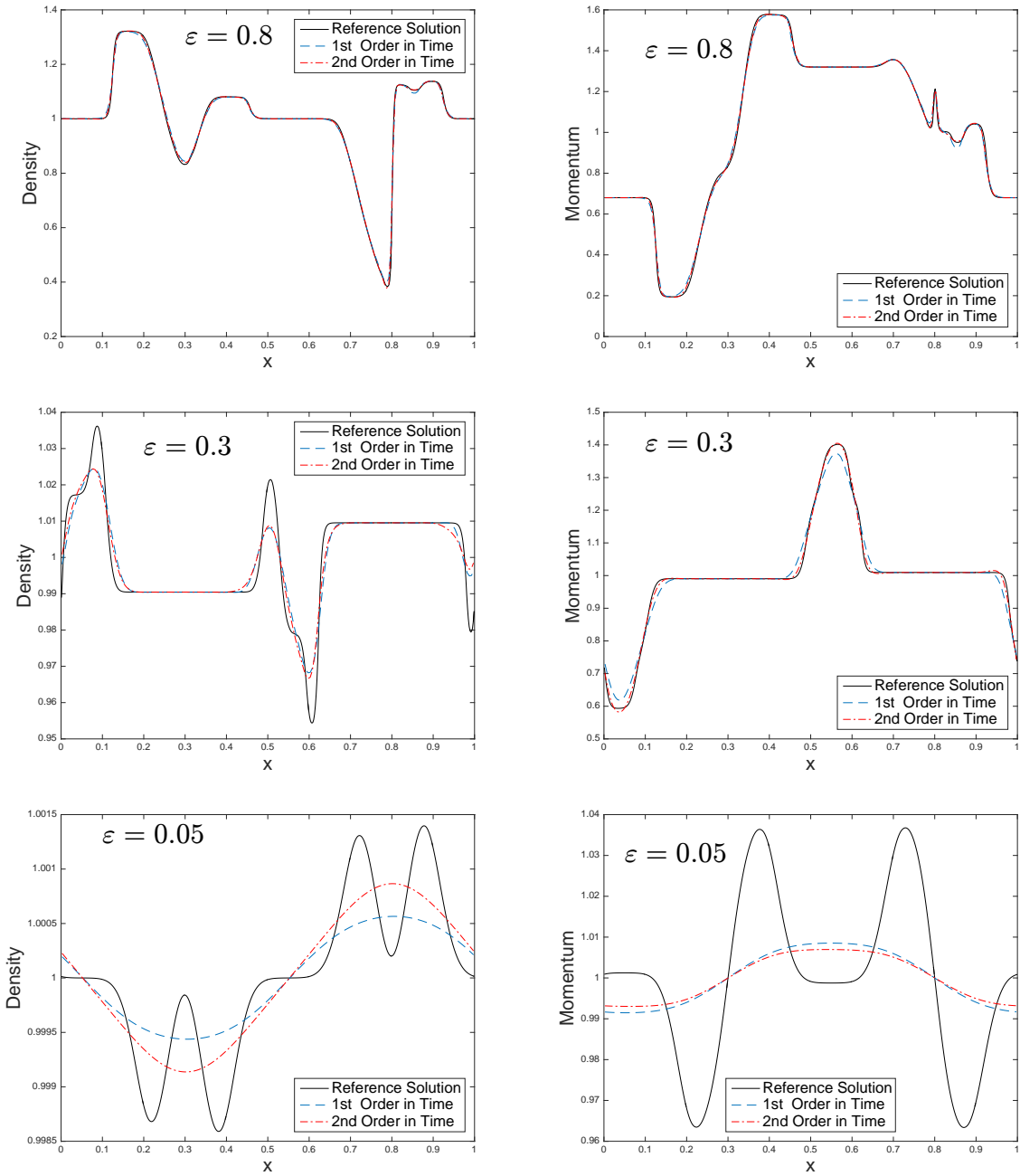

Figure 5.1: Numerical results for the Riemann problem at final time $T = 0.05$ with $\Delta x = 1/200$, $\Delta t$ is computed by (5.4) with $CFL_{Imp} = 0.5$, for the density (left) and momentum (right). The solid line is the reference solution. The corresponding values of the classical CFL numbers are: 0.3838, 0.5839, 2.9317, for $\varepsilon = 0.8$, 0.3 and 0.05, respectively.

| Density error with $CFL_{Imp} = 0.45$ and $T = 0.3$ | | | | | |
|---|---|---|---|---|---|
| | $\varepsilon = 0.8$ | | $\varepsilon = 0.3$ | | $\varepsilon = 0.05$ |
| N | $L^1$ error | $L^1$ order | $L^1$ error | $L^1$ order | $L^1$ error | $L^1$ order |
| 10 | 1.123e-02 | 0.0000 | 1.209e-02 | 0.0000 | 1.416e-04 | 0.0000 |
| 20 | 1.985e-03 | 2.4994 | 6.837e-03 | 0.8227 | 2.110e-04 | -0.5757 |
| 40 | 9.634e-04 | 1.0433 | 2.690e-03 | 1.3460 | 1.282e-03 | -2.6035 |
| 80 | 2.241e-04 | 2.1040 | 8.552e-04 | 1.6530 | 6.219e-03 | -2.2782 |
| 160 | 5.561e-05 | 2.0106 | 2.373e-04 | 1.8494 | 4.551e-03 | 0.4507 |
| 320 | 1.353e-05 | 2.0397 | 6.023e-05 | 1.9783 | 1.597e-03 | 1.5111 |
| 640 | 3.269e-06 | 2.0489 | 1.509e-05 | 1.9971 | 4.993e-04 | 1.6771 |
| 1280 | 7.868e-07 | 2.0547 | 3.774e-06 | 1.9993 | 1.358e-04 | 1.8789 |
| 2560 | 1.898e-07 | 2.0517 | 9.420e-07 | 2.0023 | 3.405e-05 | 1.9951 |

Table 5.1: Convergence table of Example 2.

denoting $U$ the numerical solution and considering the $L^1$ norm of the relative error between the numerical solution at $N$ cells and the numerical solution at $2N$ cells with

$$e_N = \frac{\left\|U_N - \overline{U}_N\right\|_{L^1}}{\left\|\overline{U}_N\right\|_{L^1}}$$

where the componentwises of the vector $\overline{U}_N$ are computed as

$$\overline{U}_N(i) = \frac{U_{2N}(2i-1) + U_{2N}(2i)}{2}.$$

The results show second-order convergence as expected for large values of $\varepsilon = 1$, and 0.3 and for small one $\varepsilon = 0.05$. We observe that because of ... of the underresolved acoustic waves, for small values of $\varepsilon$ convergence is reached for large values of $N$. We get similar results for the other variables, velocity and pressure, but we omit to show them.

## 5.3 Example 3: Two colliding acoustic waves.

Consider the evolution of two colliding acoustic waves, with the following well prepared initial data (see fig. 5.3, 5.3), i.e. when $\varepsilon$ goes to 0, the density and momentum are consistent with the incompressible limit:

$$\begin{aligned} p(\rho_\varepsilon) &= \rho_\varepsilon^\gamma, \quad \text{for} \quad x \in [-1,1], \quad \text{with} \quad \gamma = 1.4 \\ \rho_\varepsilon(x,0) &= 0.955 + \tfrac{\varepsilon}{2}(1 - \cos(2\pi x)) \quad , \quad u_\varepsilon(x,0) = -\mathrm{sign}(x)\sqrt{\gamma}(1 - \cos(2\pi x)). \end{aligned} \quad (5.6)$$

These acoustic pulses, one right-running and one left-runnig, collide and superpose and then separate again and during the whole procedure no shock forms. Now, as done in [16], we choose as spatial step $\Delta x = 1/50$ and time is computed by (5.4) with $CFL_{Imp} = 0.5$, we used $\varepsilon = 0.1$ and periodic boundary conditions. In the Figure (5.3), we display numerical results of the density and the momentum for different final times $T$. The solid line is the reference solution computed with $\Delta x = 1/500$ and $\Delta t = 1/10000$.

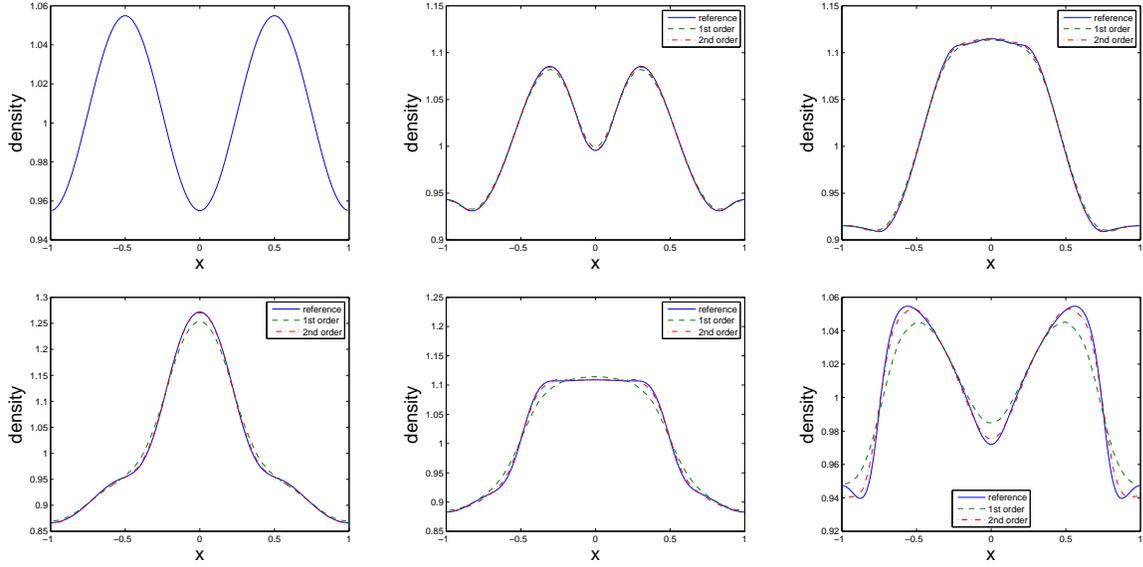

Figure 5.2: Example 3. The density computed by the first order "- -", and second order scheme "·-·" at different times: $T = 0.0$ (initial density), $T = 0.01$, $T = 0.02$, $T = 0.04$, $T = 0.06$, $T = 0.08$ with $\Delta x = 1/50$, $\varepsilon = 0.1$.

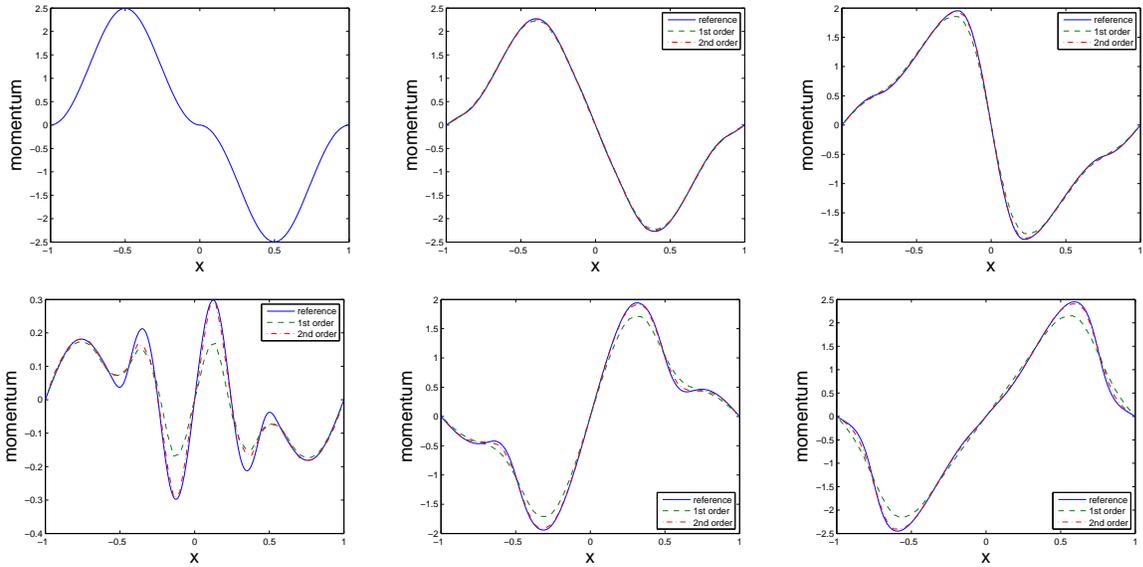

Figure 5.3: Example 3. The momentum computed by the first order "- -", and second order scheme "·-·" at different times: $T = 0.0$ (initial momentum), $T = 0.01$; $T = 0.02$; $T = 0.04$; $T = 0.06$; $T = 0.08$ with $\Delta x = 1/50$, $\varepsilon = 0.1$.

## 5.4 Example 4: 2D isentropic problem

In this section we focus on a a two dimensional case with $p(\rho) = \rho^2$. The computation domain is $\Omega = [0, 1] \times [0, 1]$ discretized by a uniform grid. We consider the same numerical test used in [19] with initial conditions:

$$\begin{cases} \rho(x,y,0) = 1 + \varepsilon^2 \sin^2(2\pi(x+y)), \\ \rho(x,y,0)u(x,y,0) = \sin(2\pi(x-y)) + \varepsilon^2 \sin(2\pi(x+y)), \\ \rho(x,y,0)v(x,y,0) = \sin(2\pi(x-y)) + \varepsilon^2 \cos(2\pi(x+y)). \end{cases}$$

In this test, the CFL condition in two dimensions is similar to (5.4). Note that at the leading order, i.e. $\varepsilon = 0$, the velocity field is divergence free and the density field is constant. In Figures 5.4 and 5.5, we display the performance of the scheme in the incompressible regime with under-resolved meshes for $\varepsilon = 0.05$ and $0.001$. We use $\Delta x = \Delta y = 1/40$, and $\Delta t$ computed by (5.4) with $CFL_{Imp} = 0.5$ at final time $T = 1$.

In Fig. 5.4 we plot the deviation from the mean of the initial density $\rho_0$ (left panels) and of the final density $\rho$ (right panels), respectively for $\varepsilon = 0.05$ and $\varepsilon = 10^{-3}$. In Fig. 5.5 we report the value of $\nabla \cdot \mathbf{u} \approx D_x u + D_y v$ at final time $T = 1$.

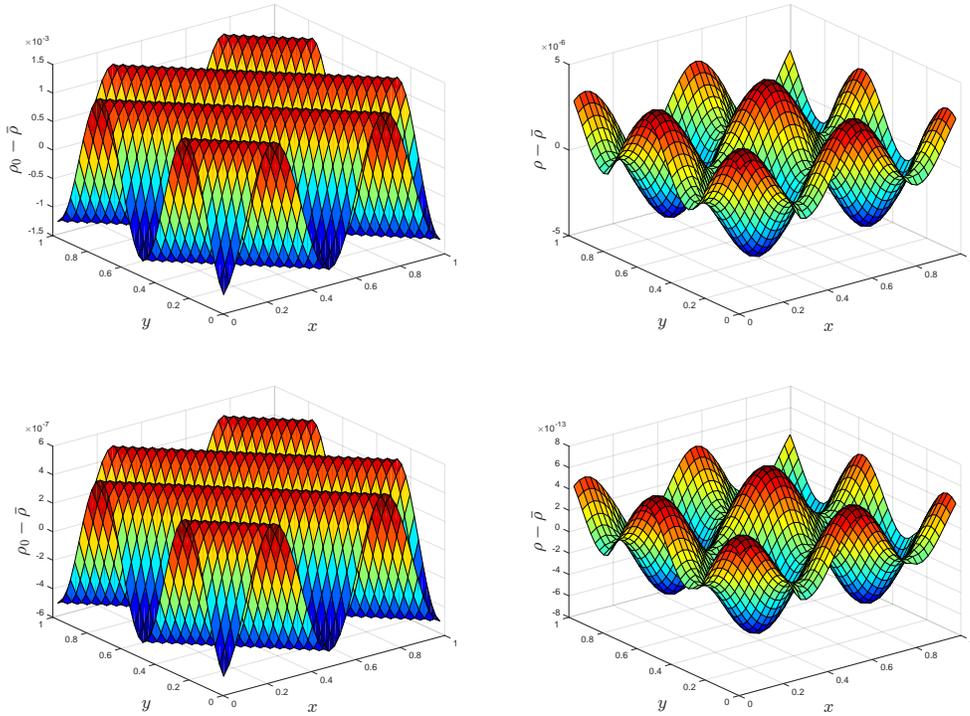

Figure 5.4: Deviation from the mean of the initial density (left panels) and of the final density at time $T = 1$ (right panels), for $\varepsilon = 0.05$ (top panels), and $\varepsilon = 10^{-3}$ (bottom panels) with $\Delta x = \Delta y = 1/40$.

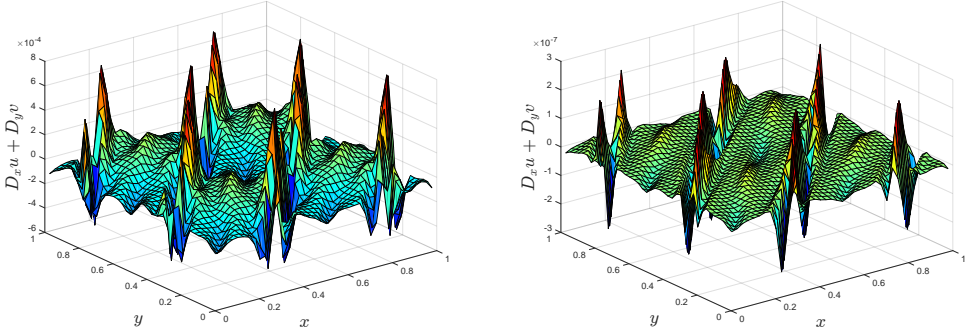

Figure 5.5: Numerical results at time $T = 1$ a): for $\varepsilon = 0.05$ (left), and $\varepsilon = 10^{-3}$ (right) with $\Delta x = \Delta y = 1/40$.

# 6 Extension to the full Euler system

In this section we consider the rescaled (non-dimensionalised) compressible Euler equations for an ideal gas (2.2) with the (suitably scaled) equation of state (EOS):

$$p = (\gamma - 1)\left(E - \varepsilon^2 \frac{\rho|\mathbf{u}|^2}{2}\right). \tag{6.1}$$

## 6.1 Reformulation of the main problem.

In order to solve numerically system (2.2), we rewrite such system in an equivalent way. We substitute the pressure (6.1) in the equation for the momentum and for the Energy $E$ and we get

$$\begin{cases} \rho_t + \nabla \cdot \mathbf{m} = 0 \\ \mathbf{m}_t + \nabla \cdot \left(\frac{\mathbf{m} \otimes \mathbf{m}}{\rho}\right) - \frac{\gamma - 1}{2}\nabla\left(\frac{|\mathbf{m}|^2}{\rho}\right) + \frac{\gamma - 1}{\varepsilon^2}\nabla E = 0 \\ E_t - \nabla \cdot \left(\frac{\gamma - 1}{2}\varepsilon^2 \frac{|\mathbf{m}|^2 \mathbf{m}}{\rho^2}\right) + \nabla \cdot (\gamma E \mathbf{u}) = 0. \end{cases} \tag{6.2}$$

Note that, the formal derivation of the incompressible Euler equations from system (6.2) is similar to the one shown in the isentropic case. We briefly summarize it here. Let us consider the following asymptotic ansatz:

$$\begin{aligned} p(\mathbf{x}, t) &= p_0(\mathbf{x}, t) + \varepsilon^2 p_2(\mathbf{x}, t) + \cdots \\ u(\mathbf{x}, t) &= u_0(\mathbf{x}, t) + \varepsilon^2 u_2(\mathbf{x}, t) + \cdots, \\ E(\mathbf{x}, t) &= E_0(\mathbf{x}, t) + \varepsilon^2 E_2(\mathbf{x}, t) + \cdots. \end{aligned} \tag{6.3}$$

Inserting the expressions into system (6.2), we get to $\mathcal{O}(\varepsilon^{-2})$, $\nabla E_0 = 0$ and this implies $E_0(x, t) = E_0(t)$, i.e. the leading order, the energy term (and hence by (6.1), the pressure $p_0 = (\gamma - 1)E_0)$ is constant in space. Furthermore, from the energy equation we get to $\mathcal{O}(1)$:

$$\nabla \cdot \mathbf{u_0} = \frac{1}{\gamma p_0} \frac{dp_0}{dt}.$$

Integrating this equation on a bounded domain $\Omega$ and using, for example, periodic boundary conditions we have: $\nabla \cdot \mathbf{u_0} = \mathbf{0}$, i.e. $p_0$ is a constant of order 1.

Furthermore, from the momentum equation we obtain, to the order $\mathcal{O}(1)$:

$$\partial_t(\rho_0 \mathbf{u_0}) + \nabla \cdot (\rho_0 \mathbf{u_0} \otimes \mathbf{u_0}) - \frac{(\gamma-1)}{2}\nabla \cdot \left(\rho_0 \mathbf{u}_0^2\right) + \frac{(\gamma-1)}{2}\nabla E_2 = 0, \qquad (6.4)$$

and by (6.1), at the $\mathcal{O}(1)$ for the pressure one has:

$$E_2 = \frac{p_2}{(\gamma-1)} + \frac{(\gamma-1)}{2}(\rho_2 \mathbf{u}_0^2). \qquad (6.5)$$

Thus, for low Mach number (i.e., $\varepsilon \ll 1$), under suitable *well- prepared* initial conditions:

$$\begin{aligned}
\rho(\mathbf{x}, t=0) &= \rho^*(x) + \varepsilon^2 \rho_2(\mathbf{x}), \\
p(\mathbf{x}, t=0) &= p_0 + \varepsilon^2 p_2(\mathbf{x}), \\
u(\mathbf{x}, t=0) &= \hat{\mathbf{u}}_0(x) + \mathcal{O}(\varepsilon),
\end{aligned} \qquad (6.6)$$

the solution $(\rho, \mathbf{m}, p)$ of (6.2), is close to the solution of the incompressible Euler system:

$$\begin{cases}
\rho_t + \mathbf{u} \cdot \nabla \rho = 0, \\
\rho(t=0, x) = \rho^*(x), \\
\partial_t(\rho \mathbf{u}) + \nabla \cdot (\rho \mathbf{u} \otimes \mathbf{u}) + \nabla p_2 = 0, \\
\nabla \cdot \mathbf{u} = 0, \\
\mathbf{u}(t=0, x) = \hat{\mathbf{u}}(x), \\
p = p_0,
\end{cases} \qquad (6.7)$$

where $\hat{\mathbf{u}}(x)$ is of order 1 with $\nabla \cdot \hat{\mathbf{u}} = 0$, $\rho^*(x)$ is a strictly positive function such that $\rho^* = \mathcal{O}(1)$, and $p_0$ is the space independent thermodynamic pressure, related to internal energy $e_0$ and density $\rho_0$ by: $p_0 = (\gamma-1)\rho_0 e_0$. Note that $p_2 = \lim_{\varepsilon \to 0} \frac{1}{\varepsilon^2}(p-p_*)$ is implicitly defined by the constraint $\nabla \cdot \mathbf{u} = 0$ and explicitly given by the equation $-\Delta p_2 = \rho_0 \nabla^2 : (u \otimes u)$.

## 6.2 The semi-implicit R-K (SI R-K) time integrator for the Euler equation

The scheme proposed here for the full Euler system follows the isentropic case strategy. We start applying to system (6.2) an explicit-implicit Euler scheme:

$$\begin{cases} \dfrac{\rho^{n+1} - \rho^n}{\Delta t} + \nabla \cdot \mathbf{m}^{n+1} = 0, \\ \dfrac{\mathbf{m}^{n+1} - \mathbf{m}^n}{\Delta t} + \nabla \cdot \left( \dfrac{\mathbf{m}^n \otimes \mathbf{m}^n}{\rho^n} \right) - \dfrac{\gamma - 1}{2} \nabla \left( \dfrac{(|\mathbf{m}|^2)^n}{\rho^n} \right) + \dfrac{\gamma - 1}{\varepsilon^2} \nabla E^{n+1} = 0, \\ \dfrac{E^{n+1} - E^n}{\Delta t} - \nabla \cdot \left( \dfrac{\gamma - 1}{2} \varepsilon^2 \dfrac{(|\mathbf{m}|^2 \mathbf{m})^n}{\rho^n} \right) + \gamma \nabla \cdot \left( E^n \dfrac{\mathbf{m}^{m+1}}{\rho^n} \right) = 0. \end{cases} \quad (6.8)$$

with

$$p^{n+1} = (\gamma - 1) \left( E^{n+1} - \varepsilon^2 \dfrac{(\rho \mathbf{u}^2)^n}{2} \right). \quad (6.9)$$

Then, the scheme works as follows. Let

$$\rho^{n+1} = \rho^n - \Delta t \nabla \cdot \mathbf{m}^{n+1},$$

$$\overset{1}{m}{}^{n+1} = \overset{1}{m}{}^* - \dfrac{(\gamma - 1)}{\varepsilon^2} \Delta t \partial_x E^{n+1},$$

$$\overset{2}{m}{}^{n+1} = \overset{2}{m}{}^* - \dfrac{(\gamma - 1)}{\varepsilon^2} \Delta t \partial_y E^{n+1}, \quad (6.10)$$

$$E^{n+1} = E^* - \gamma \Delta t \partial_k \left( E^n \dfrac{\overset{k}{m}{}^{n+1}}{\rho^n} \right),$$

with

$$\overset{k}{m}{}^* = \overset{k}{m}{}^n - \Delta t \nabla \cdot \left( \dfrac{\overset{k}{m}{}^n \mathbf{m}^n}{\rho^n} \right) + \dfrac{\gamma - 1}{2} \Delta t \partial_k \left( \dfrac{(|m|^2)^n}{\rho^n} \right), \quad k = 1, 2, \quad (6.11)$$

and

$$E^* = E^n - \dfrac{\gamma - 1}{2} \varepsilon^2 \Delta t \partial_k \left( \dfrac{|m^2|^n \overset{k}{m}{}^n}{\rho^n} \right), \quad k = 1, 2, \quad (6.12)$$

Now, inserting the expression $\overset{k}{m}{}^{n+1}$ into the equation for $E$ we obtain

$$E^{n+1} = E^{**} + \dfrac{\Delta t^2 \gamma (\gamma - 1)}{\varepsilon^2} \nabla \cdot \left( \dfrac{E^n}{\rho^n} \nabla \cdot E^{n+1} \right), \quad (6.13)$$

where

$$E^{**} = E^* - \Delta t \gamma \partial_k \left( \dfrac{E^n}{\rho^n} \overset{k}{m}{}^* \right). \quad (6.14)$$

Then discretizing again in space by central scheme based on a staggered grid with uniform spatial mesh, we obtain in this case a linear elliptic equation:

$$E_{i+1/2,j+1/2}^{n+1} = E_{i+1/2,j+1/2}^{**} + \frac{\Delta t^2}{\varepsilon^2}\gamma(\gamma-1)LE_{i+1/2,j+1/2}^{n+1}, \tag{6.15}$$

Here the operator $L$ is defined as: $LE_{ij} = L_x E_{ij} + L_y E_{ij}$ where:

$$L_x E_{ij}^{n+1} = \frac{1}{\Delta x}\left(\left(\frac{E}{\rho}\right)_{i+j/2,j}^n \left(\frac{E_{i+1,j}^{n+1} - E_{i,j}^{n+1}}{\Delta x}\right) - \left(\frac{E}{\rho}\right)_{i-j/2,j}^n \left(\frac{E_{i,j}^{n+1} - E_{i-1,j}^{n+1}}{\Delta x}\right)\right),$$

where for any function $k$ defined on the grid:

$$k_{i+1/2,j} = \frac{k_{i+1,j} + k_{i,j}}{2},$$

and similarly for $L_y$. Then, solving the linear system, we obtain the new energy $E$ to step $n+1$. This is used to update the momentum $m^{n+1}$, and, finally, to compute the new density.

**AP property.** Following the analysis in Section 3.2, it is possible to show that scheme (6.8) is asymptotic preserving (AP) in the limit $\varepsilon \to 0$. We start by making the following ansatz: the symptotic *ansatz*

$$\begin{aligned} p_{ij}^n &= p_0^n + \varepsilon^2 p_{2,ij}^n + \cdots, \\ \rho_{ij}^n &= \rho_{0,ij}^n + \varepsilon^2 \rho_{2,ij}^n + \cdots, \\ E_{ij}^n &= E_0^n + \varepsilon^2 E_{2,ij}^n + \cdots. \end{aligned} \tag{6.16}$$

with the leading order energy term: $E_0^0 = p_0^n/(\gamma-1)$ and from (6.5) and (6.9) we get to $\mathcal{O}(1)$:

$$E_{2,ij}^n = \left(\frac{p_{2,ij}^n}{(\gamma-1)} + \frac{\rho_{0,ij}^{n-1}|\mathbf{u}_{0,ij}^{n-1}|^2}{2}\right). \tag{6.17}$$

Furthermore, in agreement with (6.6), we consider the following *well-prepared* initial conditions

$$\begin{aligned} \rho^0 &= \rho_{0,ij}^* + \varepsilon^2 \rho_{2,ij}^0, \\ p^0 &= p_0 + \varepsilon^2 p_{2,ij}^0, \\ u^0 &= u_{0,ij}^0 + \mathcal{O}(\varepsilon), \end{aligned} \tag{6.18}$$

with $p_0$ a positive constant, and $\rho_{0,ij}^*$, $u_{0,ij}^0$ functions independent on time. Then following the same procedure adopted in Sect. 3.2 one obtains a discretization of sytem (6.7). The reader can find a detailed asymptotic preserving analysis of the scheme (6.8) in [44].

**High order extension.** A natural extension of the first order scheme (6.8) to a second order one is obtained in a similar way to the one in Sect. 4. In particular we obtain a second order globally stiff accurate SI-IMEX-RK-STAG. Here the function $\mathcal{H}(\mathbf{U}^*, \mathbf{U})$ is defined, by system (6.2), as follows:

$$\mathcal{H}(\mathbf{U}^*, \mathbf{U}) = \begin{pmatrix} -\nabla \cdot \mathbf{m} \\ -\nabla \cdot \left(\frac{\mathbf{m}^* \otimes \mathbf{m}^*}{\rho^*}\right) + \frac{\gamma-1}{2}\nabla\left(\frac{|\mathbf{m}^*|^2}{\rho^*}\right) - \frac{\gamma-1}{\varepsilon^2}\nabla E \\ \nabla \cdot \left(\frac{\gamma-1}{2}\varepsilon^2 \frac{|\mathbf{m}^*|^2 \mathbf{m}^*}{\rho^*}\right) - \gamma\nabla \cdot \left(\frac{E^*}{\rho^*}\mathbf{m}\right). \end{pmatrix} \tag{6.19}$$

# 7 Numerical Results for the full Euler system.

In this last part of the paper we present several numerical test cases applied to the full Euler system in one and two dimensions. Here we consider the same CFL condition and time step given in Sect. 5.1.

## 7.1 Example 1: 1D Sod Shock Tube Problem.

First, we test the method in a compressible regime, i.e. when the Mach number is $\mathcal{O}(1)$. We consider the classical Sod shock tube problem with the initial conditions

$$(\rho, u, p)(x, 0) = \begin{cases} (1.0,\ 0.0,\ 1.0) & \text{if} \quad x < 0.5, \\ (0.125,\ 0.0,\ 0.1) & \text{otherwise.} \end{cases} \tag{7.1}$$

The domain is $[0, 1]$ and the discontinuity is initially at $x = 0.5$ we consider $N = 200$ grid points. The numerical results are obtained by the second order scheme at $T = 0.18$ and are shown in Figure 7.1. The reference solution is computed with $N = 1000$ grid points.

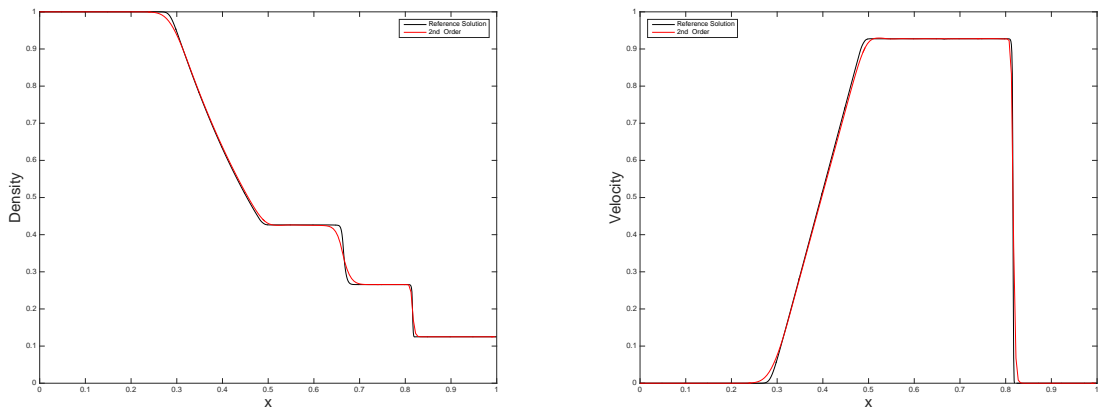

Figure 7.1: Sod shock tube, Density and velocity solution (red lines), $CFL = 0.5$, $N = 200$ grid points at time $T = 0.18$ with $\varepsilon = 1$. Black line is the reference solution.

The results are comparable to those obtained by Nessyahu-Tadmor central scheme, [36].

## 7.2 Example 2: Convergence test.

We compute the experimental order of convergence (EOC) of the scheme presented in the previous section. The initial conditions, the final time and the boundary conditions are given in Sect. 5.2. Here we use the value $\gamma = 1.4$.

Below we report the convergence table for the density, momentum and energy for different values of the Mach number 0.8, 0.3 and $10^{-4}$. Note that for $\varepsilon = 10^{-4}$ we choose as final time $T = 0.01$.

| $CFL_{Imp} = 0.45$, $T = 0.3$ and $\varepsilon = 0.8$ | | | | | |
|---|---|---|---|---|---|
| N | $L^1$-error $\rho$ | EOC $\rho$ | $L^1$-error $m$ | EOC $m$ | $L^1$-error $E$ | EOC $E$ |
| 20 | 5.472e-03 | 0.0000 | 1.347e-02 | 0.0000 | 6.990e-03 | |
| 40 | 1.602e-03 | 1.7727 | 4.047e-03 | 1.7352 | 2.603e-03 | 1.4251 |
| 80 | 4.792e-04 | 1.7408 | 1.302e-03 | 1.6361 | 7.481e-04 | 1.7990 |
| 160 | 1.237e-04 | 1.9536 | 3.639e-04 | 1.8392 | 2.018e-04 | 1.8900 |
| 320 | 3.120e-05 | 1.9874 | 9.381e-05 | 1.9556 | 5.146e-05 | 1.9715 |
| 640 | 7.732e-06 | 2.0126 | 2.353e-05 | 1.9955 | 1.286e-05 | 2.0002 |

| $CFL_{Imp} = 0.45$, $T = 0.3$ and $\varepsilon = 0.1$ | | | | | |
|---|---|---|---|---|---|
| N | $L^1$-error $\rho$ | EOC $\rho$ | $L^1$-error $m$ | EOC $m$ | $L^1$-error $E$ | EOC $E$ |
| 20 | 1.988e-03 | 0.0000 | 3.894e-02 | 0.0000 | 2.708e-03 | |
| 40 | 3.959e-04 | 2.3281 | 7.445e-03 | 2.3870 | 5.501e-04 | 2.2998 |
| 80 | 1.202e-04 | 1.7199 | 2.241e-03 | 1.7319 | 1.676e-04 | 1.7145 |
| 160 | 3.036e-05 | 1.9850 | 5.641e-04 | 1.9904 | 4.266e-05 | 1.9742 |
| 320 | 7.628e-06 | 1.9929 | 1.400e-04 | 2.0100 | 1.077e-05 | 1.9863 |
| 640 | 1.895e-06 | 2.0090 | 3.501e-05 | 2.0000 | 2.682e-06 | 2.0051 |

| $CFL_{Imp} = 0.45$, $T = 0.01$ and $\varepsilon = 10^{-4}$ | | | | | |
|---|---|---|---|---|---|
| N | $L^1$-error $\rho$ | EOC $\rho$ | $L^1$-error $m$ | EOC $m$ | $L^1$-error $E$ | EOC $E$ |
| 20 | 1.844e-05 | 0.0000 | 4.792e-01 | 0.0000 | 2.582e-05 | |
| 40 | 1.084e-05 | 0.7667 | 2.105e-01 | 1.1865 | 1.518e-05 | 0.7667 |
| 80 | 2.929e-06 | 1.8881 | 5.463e-02 | 1.9463 | 4.100e-06 | 1.8881 |
| 160 | 7.332e-07 | 1.9980 | 1.360e-02 | 2.0057 | 1.026e-06 | 1.9980 |
| 320 | 1.831e-07 | 2.0018 | 3.394e-03 | 2.0031 | 2.563e-07 | 2.0018 |
| 640 | 4.582e-08 | 1.9982 | 8.492e-04 | 1.9986 | 6.415e-08 | 1.9982 |

We observe second order convergence in all case. However, because acoustic waves are poorly revolved in time, we need to use a very fine mesh in order to observe expected order of convergence.

### 7.3 Example 2: Two Colliding Acoustic Pulses.

We consider again two colliding acoustic pulses in a weakly compressible regime. This test has been taken from [37, 31]. The initial conditions are given by

$$\begin{aligned}
p(\rho_\varepsilon) &= \rho_\varepsilon^\gamma, \quad \text{for} \quad x \in [-1, 1], \quad \text{with} \quad \gamma = 1.4 \\
\rho_\varepsilon(x, 0) &= \rho_0 + \tfrac{1}{2}\varepsilon\rho_1(1 - \cos(2\pi x)) \quad \rho_0 = 0.955, \quad \rho_1 = 2.0 \\
u_\varepsilon(x, 0) &= \tfrac{1}{2}u_0\text{sign}(x)u_0(1 - \cos(2\pi x)) \quad u_0 = 2\sqrt{\gamma} \\
p(x, 0) &= p_0 + \tfrac{1}{2}\epsilon p_1(1 - \cos(2\pi x)), \quad p_0 = 1.0, \quad p_1 = 2\gamma
\end{aligned} \quad (7.2)$$

The domain is $-L \leq x \leq L = 2/\varepsilon$ and the boundary conditions are periodic. In Figure 7.3 we show the plots of the pressure obtained using first order scheme (6.8) and the second order globally stiffly accurate SI-IMEX-RK-STAG scheme at different times $t = 0.815$ and $t = 1.63$. We choose $\varepsilon = 1/11$ and we plot the initial pressure distributions for comparison. In the figure, the reference solution has been computed with $N = 1500$. The results have the same behaviour as in [37].

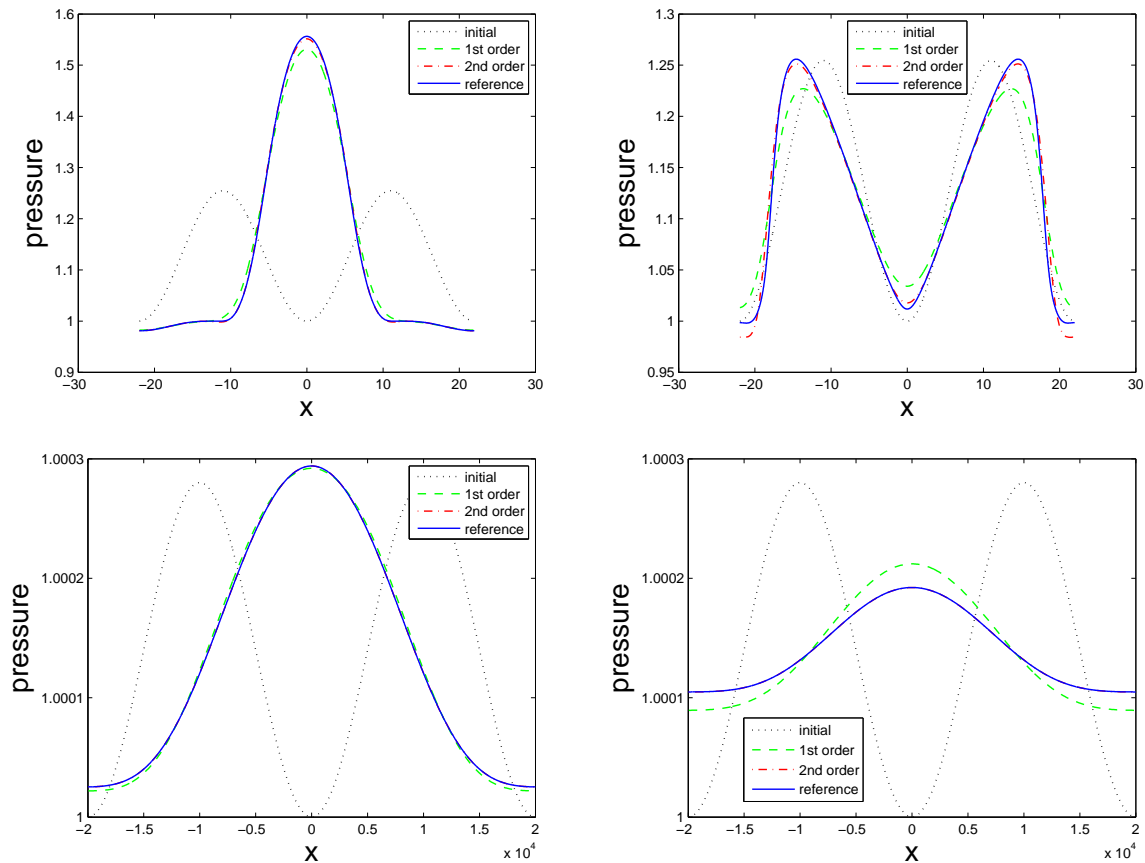

Figure 7.2: $\varepsilon = 1/11$ (top panel) and $\varepsilon = 10^{-4}$ (bottom panel), $N = 440$, final time $T = 0.815$ (Left) and $T = 1.63$ (Right). We displayed the initial condition "$\cdot \cdot \cdot$", the numerical solution to the first "- -" and second order scheme "$\cdot$ - $\cdot$". Reference solution "__".

## 7.4 Example 3: Asymptotic Preserving Property (AP).

In Section 3.2 we have discussed the asymptotic preserving property of our scheme when applied to the Euler isentropic case (2.3) and to the full Euler one (6.2). Now, the goal in this section is to verify such AP property, numerically.

In particular, we compare the solutions of the scheme at low Mach number, when applied to the isentropic Euler equations (2.3) and to the complete Euler ones (6.2), with the solution of the incompressible Euler equations.

Consider incompressible Euler equation in the vorticity stream-function formulation, [43]:

$$\omega_t + \mathbf{u} \cdot \nabla \omega = 0, \quad (x,y) \in [0, 2\pi] \times [0, 2\pi] \tag{7.3}$$

where

$$\omega = \frac{\partial v}{\partial x} - \frac{\partial u}{\partial y}.$$

Because $\nabla \cdot \mathbf{u} = 0$, there exists a function $\psi$ such that $\mathbf{u} = (\partial_y \psi, -\partial_x \psi)$. Inserting this relation in the expression for $\omega$ one obtains the Poisson equation $-\Delta \psi = \omega$.

For our numerical test we assume periodic boundary conditions and we consider as initial condition, the shear flow:

$$\omega(x,y,0) = \begin{cases} \delta \cos(x) - \frac{1}{\rho}\text{sech}^2((y-\pi/2)/\rho), & y \leqslant \pi \\ \delta \cos(x) + \frac{1}{\rho}\text{sech}^2((3\pi/2-y)/\rho), & y > \pi \end{cases} \tag{7.4}$$

where $\delta = 0.05$ and $\rho = \pi/15$.

As a reference solution for incompressible Euler equations (7.3), we consider a very accurate numerical solution obtained by Fourier spectral discretization in space and fourth order Runge-Kutta method in time. Final time is $T = 6$ with $N = 160$ grid points we have a fully resolved calculation. We refer it as *reference* solution (see Fig. 7.3).

Now we compute numerical solutions for isentropic Euler equations (2.3) and to the complete Euler ones (6.2) by using our scheme. We set $\varepsilon = 10^{-4}$, initial conditions (7.4), periodic boundary conditions, $N = 160$ and final time $T = 6.0$, the results are given in Figures 7.4.

Note that the CFL condition in two dimensions is similar to (5.4), i.e.,

$$\Delta t = \frac{CFL_{Imp}}{\frac{\max |\mathbf{u}_{ij} + \tilde{c}_{ij}|}{\Delta x} + \frac{\max |\mathbf{v}_{ij} + \tilde{c}_{ij}|}{\Delta y}}.$$

In our simulation we used $\Delta x \equiv \Delta y$.

Comaring the reference solution for incompressible Euler equations, Fig. 7.3, with the solutions obtained with the second order semi implicit scheme (SI-IMEX-RK-STAG), Fig. 7.4, we observe that there is a qualitative agreement.

Finally, in Figure 7.5 we show the behaviour of the $L^1$ error, as the difference between the numerical solution of the Euler scheme computed by our second order scheme, SI-IMEX-RK-STAG, with different values of $\varepsilon$ from $10^{-1}$ to $10^{-4}$ and the reference solution. We refined $N = 16, 32, 64, \cdots, 1024$ to compute the numerical solution for each value of $\varepsilon$, instead we fixed $N = 1024$ for the reference one. The final time is $T = 1$.

Two main sources of error are presented when comparing the numerical solution of weakly compressible Euler equation and the reference solution for incompressible Euler: space and time discretization error and modelling error. The latter is due to the different behaviour between compressible and incompressible flow. For a given value of $\varepsilon$, discretization error is dominant for coarse grid. Refining the grid the typical second order convergence rate is observed. When space and time are fully resolved, the model error becomes dominant and is responsible of the plateau observed for large values of $N$. As $\varepsilon$ is decreased, the discrepancy between compressible and incompressible Euler decreases, and is observed at finer and finer meshes (see Fig. 7.5).

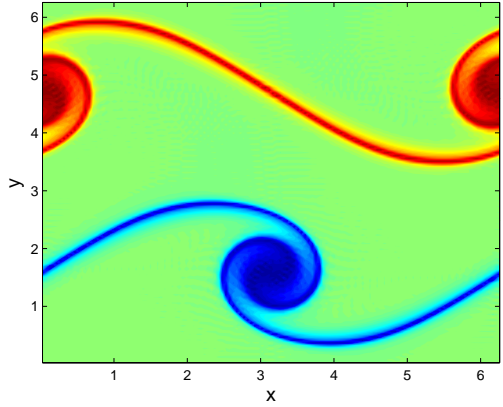

Figure 7.3: Plot of the reference solution at final time $T = 6$.

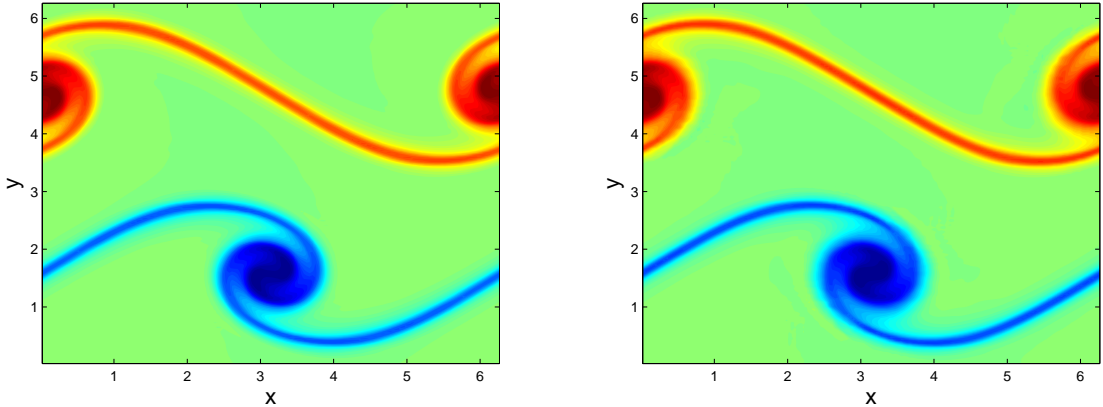

Figure 7.4: Numerical solution for the isentropic case (left), and for the complete Euler case (right) with $\varepsilon = 10^{-4}$ at final time $T = 6$.

# 8 Conclusions

In this paper we present a simple strategy to construct semi-implicit schemes for Euler equations which are able to work on a wide range of Mach numbers. The schemes are based on staggered grid discretization in space, which is second order accurate and avoids the needs of exact or approximate Riemann solvers. In the proposed schemes, acoustic waves are treated linearly implicitly, while material waves are treated explicitly. The resulting scheme is Asymptotic Preserving, i.e. it converges to a consistent scheme for the incompressible Euler equations as the Mach number vanishes. Because of the linear treatment of implicit terms, the schemes are quite efficient, especially for low Mach number flow. Second order schemes in time are constructed using Implicit-Explicit Runge-Kutta methods. Because of the staggered nature of the problem, second order discretization in time requires the use of Globally Stiffly Accurate schemes. A simplification is obtained by computing most of the stage values implicitly by a possibly non conservative predictor, while the conservative corrector guarantees that the overall scheme is conservative. A stability analysis is performed showing that the last stage of the scheme has to be implicit at the level of the numerical solution in order to avoid classical stability restrictions. Numerical convergence study is performed on various test problems, emphasizing the robustness and efficiency of the scheme. The procedure could be extended to construct higher order schemes. This issue, together with implicit treatment of boundary conditions, will be subject of future investigation.

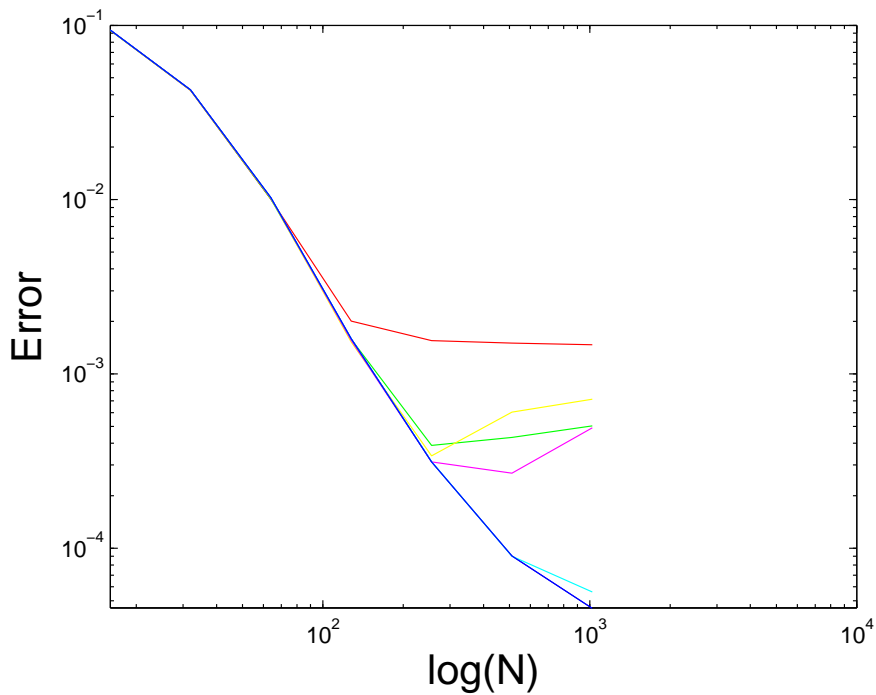

Figure 7.5: $L^1$ error for: $\varepsilon = 10^{-1},\ 7.7 \cdot 10^{-2},\ 5.5 \cdot 10^2,\ 3.25 \cdot 10^{-2}, 10^{-3},\ 10^{-4}$ and $N = 16, 32, 64, \cdots, 1024$. The final time is $T = 1$.

# Acknowledgments

The work has been partially supported by ITN-ETN Horizon 2020 Project *ModComp-Shock, Modeling and Computation on Shocks and Interfaces*, Project Reference 642768,